\documentclass[12pt]{article}

\usepackage{times}
\usepackage[centertags]{amsmath}
\usepackage{amssymb}
\usepackage{amsthm}
\usepackage[latin1]{inputenc}
\usepackage[colorlinks=true,urlcolor=red,citecolor=red,linkcolor=red]{hyperref}
\usepackage[a4paper, top=5cm, bottom=5cm, left=3.5cm, right=2.5cm]{geometry}
\usepackage[italian,english]{babel}
\usepackage{graphicx}
\usepackage{tikz}
\usepackage{qtree}
\usepackage{tikz-qtree}
\usepackage[lined,algonl,boxed]{algorithm2e}
\usetikzlibrary{trees,positioning,shapes} 
\newtheorem{prop}{Proposition}[section]
\newtheorem{theo}[prop]{Theorem}
\newtheorem{theo2}[prop]{Theorem}
\newtheorem{theo3}[prop]{Theorem}
\newtheorem{question}[prop]{Question}
\newtheorem{Lemma1}[prop]{Lemma}

\newtheorem{Lemma2}[prop]{Lemma}
\newtheorem{prop2}[prop]{Proposition}
\theoremstyle{definition} 
\newtheorem{example1}[prop]{Example}
\newtheorem{example2}[prop]{Example}
\newtheorem{example23}[prop]{Example}
\newtheorem{example3}[prop]{Example}
\newtheorem{example4}[prop]{Example}
\newtheorem{example5}[prop]{Example}
\newtheorem{example6}[prop]{Example} 
\newtheorem{example67}[prop]{Example}
\newtheorem{example7}[prop]{Example}
\newtheorem{example8}[prop]{Example}
\newtheorem{example9}[prop]{Example}
\newtheorem{rem}[prop]{Remark}
\newtheorem{rem1}[prop]{Remark}
\newtheorem{remark}[prop]{Remark}
\newtheorem{remark1}[prop]{Remark}
\newtheorem{remark2}[prop]{Remark}
\newtheorem{remark23}[prop]{Remark}
\newtheorem{remark24}[prop]{Remark}

\newtheorem{remark4}[prop]{Remark}

\newenvironment{acknowledgements}%
    {\null\vfill\begin{center}%
    \bfseries Acknowledgements\end{center}}%
    {\vfill\null}

\newcounter{lastnote}

\title{Arf good semigroups}

\author
{Giuseppe Zito\\
\\
\\
}

\date{}

\begin{document}

\baselineskip14.2806pt

\maketitle

\begin{abstract}
In this paper we study the property of the Arf good subsemigroups of $\mathbb{N}^n$, with $n\geq2$. We give a way to compute all the Arf semigroups with a given collection of multiplicity branches. We also deal with the problem of determining the Arf closure of a set of vectors and of a good semigroup, extending the concept of characters of an Arf numerical semigroup  to Arf good semigroups.
\end{abstract}

\let\thefootnote\relax

\section*{Introduction}

In this paper we study a particular class of good subsemigroups of $\mathbb{N}^n$. The concept of good semigroup was introduced in \cite{BDF}. Its definition depends on the properties of the value semigroups of one dimensional analytically unramfied ring (for example the local rings of an algebraic curve), but in the same paper it is shown that the class of good semigroups is  bigger than the class of value semigroups. Therefore the good semigroups can be seen as a natural generalization of the numerical semigroup and can be studied without referring to the ring theory context, with a more combinatorical approach. In this paper we deal only with local good semigroups, i.e  good semigroups $S \subseteq \mathbb{N}^n$ such that the only element of $S$ with zero component is the zero vector.

In this paper we focus on the class of local Arf good semigroups. This is motivated by the importance of the Arf numerical semigroups in the study of the equivalence of two algebroid branches. Given an algebroid branch $R$, its multiplicity sequence is defined to be the sequence of the multiplicities of the succesive blowups $R_i$ of $R$. Two algebroid branches are equivalent if and only if they have the same multiplicity sequence (cf. \cite[Definition 1.5.11]{Campi}). In \cite{Arf} it is introduced the concept of Arf ring and it is shown that for each ring $R$ there is a smallest Arf overring $R'$, called the Arf closure of $R$, that has also the same multiplicity sequence of $R$. In the same paper it is proved that two algebroid branches are equivalent if and only if their Arf closure have the same value semigroup, that is a numerical Arf semigroup, i.e. a numerical semigroup $S$ such that $S(s)-s$ is a semigroup, for each $s\in S$, where $S(s)=\left\{n \in S; n\geq s \right\}$. All these facts can be generalized to algebroid curves (with more than one branch) and this naturally leads to define the Arf good semigroups of  $\mathbb{N}^n$ by extending the numerical definition considering the usual partial ordering given by the components. 

 In the numerical case an Arf semigroup $S=\left\{ s_0=0<s_1<s_2,\ldots \right\}$ is completely described by its multiplicity sequence, that is the sequence of the differences $s_{i+1}-s_i$. Extending the concept of multiplicity sequence, in \cite{BDF} it is also shown that to each local Arf good semigroup can be associated a multiplicity tree that characterizes the semigroup completely. A tree $T$ of vectors of $\mathbb{N}^n$ has to satisfy some properties to be a multiplicity tree of a local Arf good semigroup. For instance it must have multiplicity sequences along its branches (since the projections are Arf numerical semigroups) and each node must be able to be expressed as a sum of nodes in a subtree of $T$ rooted in it.
Thus, taking in account this $1$-$1$ correspondence, the aim of this paper is to study  Arf good semigroups by characterizing their multiplicity trees, finding an unambiguous way to describe them.
Using this approach, we can also deal with the  problem of finding the Arf closure of a good semigroup $S$, that is the smallest Arf semigroup containing $S$.

The structure of the paper is the following.

In Section \ref{section2}, given a collection of $n$ multiplicity sequences $E$, we define the set $\sigma(E)$ of all the Arf semigroups $S$ such that the $i$-th projection $S_i$ is an Arf numerical semigroup associated to the $i$-th multiplicity sequences of $E$. We define also the set $\tau(E)$ of the corresponding multiplicity trees and we describe a tree in $\tau(E)$ by an upper triangular matrix $\left(p_{i,j}\right)$, where  $p_{i,j}$ is the highest level where the $i$-th and $j$-th branches are glued, and we give a way to deduce from $E$ the maximal value that can be assigned to the $p_{i,j}$.
This fact let us to understand when the set $\sigma(E)$ is finite. 
We introduce the class of untwisted trees that are easier to study because they are completely described by the second diagonal of their matrix, and we notice that a tree can be always transformed in to an untwisted one by permuting its branches.

In Section \ref{section3} we address the problem of understanding when a set of vectors $G \subseteq \mathbb{N}^n$ determines uniquely an Arf semigroup of $\mathbb{N}^n$. Thus we define $\textrm{Arf}(G)$ as the minimum of the set $S(G)=\left\{ S: S\subseteq \mathbb{N}^n \textrm{ is an Arf semigroup  and } G\subseteq S \right\}$ , and we find the properties that $G$ has to satisfy in order to have a good definition for $\textrm{Arf}(G)$ (cf. Theorem \ref{T45}). Finally, given a $G$ satisfying these properties, we give a procedure for computing $\textrm{Arf}(G)$.

In Section \ref{section4} we adapt the techniques learned in the previous section to the problem of determining the Arf closure of a good semigroup. In \cite{DPMT}, the authors solved this problem for $n=2$, leaving it open for larger dimensions. In this section we use the fact that a good semigroup $S$ can be  completely described by its finite subset $ \textrm{Small}(S)=\left\{ s\in S : s\leq \delta \right\}$, where $\delta$ is the smallest element such that $\delta+\mathbb{N}^n \subseteq S$, whose existence is guaranteed by the properties of the good semigroups. 

Finally, in Section \ref{section5}, we address the inverse problem: given an Arf semigroup $S\subseteq \mathbb{N}^n$, find a set of vectors $G \subseteq \mathbb{N}^n$ , called set of generators of $S$, such that $\textrm{Arf}(G)=S,$ in order to find a possible generalization of the concept of characters in the numerical case. In Theorem \ref{T5}, we find the properties that such a $G$ has to satisfy and we focus on the problem of finding a minimal one. From this point of view we are able to give a lower and an upper bound for the minimal cardinality   for a set of generators of a given Arf semigroup (Corollary \ref{prop2}). With an example we also show that, given an Arf semigroup $S$, it is possible to  find minimal sets of generators with distinct cardinalities.

The procedures presented here have been implemented in GAP (\cite{gap}).

\section{Arf semigroups with a given collection of multiplicity branches} \label{section2}
In this section we determine all the local Arf good semigroups having the same collection of multiplicity branches. \\
First of all we need to fix some notations and recall the most important definitions. In the following, given a vector $ v$ in $\mathbb{N}^n$, we will always denote by $v[i]$ its $i$-th component.

 A good semigroup $S$ of $\mathbb{N}^n$ is a submonoid of $\left( \mathbb{N}^n,+\right)$ such that:  (cf. \cite{BDF})\begin{itemize}
\item for all $a,b \in S$, $	\min(a,b) \in S$; 
\item if $a,b\in S$ and $a[i]=b[i]$ for some $i \in \left\{ 1, \ldots,n\right\}$, then there exists $c \in S$ such that $c[i] > a[i]=b[i]$, $c[j]\geq \min(a[j],b[j])$ for $j \in  \left\{ 1, \ldots,n\right\} \setminus \left\{ i\right\}$ and $c[j]=\min(a[j],b[j])$ if $a[j]\neq b[j]$;
\item there exists $\delta \in S$ such that $\delta+\mathbb{N}^n \subseteq S$
\end{itemize}
(where we are considering the usual partial ordering in $\mathbb{N}^n$: $a \leq b$ if $a[i] \leq b[i]$ for each $i=1,\ldots,n$).
\\
In this paper we will always deal with local good semigroups. A good semigroup $S$  is local if the zero vector is the only vector of $S$ with some component equal to zero.
However, it can be shown that every good semigroup is the direct product of local semigroups (cf. \cite[Theorem 2.5]{BDF}).

An Arf semigroup of $\mathbb{N}^n$, is a good semigroup such that $ S(\alpha)-\alpha $ is a semigroup for each $\alpha \in S$  where $ S(\alpha)=\left\{\beta \in S; \beta \geq \alpha \right\}$.
The  multiplicity tree $T$ of a local Arf semigroup $S\subseteq \mathbb{N}^n$ is a tree where the nodes are vector $ \textbf{n}_i^j \in \mathbb{N}^n$ (where with $\textbf{n}_i^j$ we mean that  this node is in the $i$-th branch on the $j$-th level. The root of the tree is $\textbf{n}_{1}^1=\textbf{n}_i^1 $ for all $i$ because we are in the local case and at level one all the branches must be glued) and we have 

$$ S=\left\{\textbf{0}\right\} \bigcup_{T'} \left\{ \sum_{\textbf{n}_i^j \in T' } {\textbf{n}_i^j}\right\},$$
where $T'$ ranges over all finite subtree of $T$ rooted in   $\textbf{n}_1^1$.

Furthermore a tree $T$ is a multiplicity tree of an Arf semigroup if and only if its  nodes  satisfy the following properties (cf. \cite[Theorem 5.11]{BDF}):

\begin{itemize}
\item there exists $L \in \mathbb{N}$ such that for $m \geq L$, $\textbf{n}_i^m=(0,\ldots,0,1,0\ldots,0)$ (the nonzero coordinate is in the $i$-th position) for any $i=1,\ldots,n$;
\item $\textbf{n}_i^j[h]=0$ if and only if $\textbf{n}_i^j$ is not in the $h$-th branch of the tree;
\item each $\textbf{n}_i^j$ can be obtained as a sum of nodes in a finite subtree $T'$ of $T$ rooted in $\textbf{n}_i^j$.
\end{itemize}
Notice that from these properties it follows that we must have multiplicity sequences along each branch.
\\ \\
Suppose now that  $E$ is an ordered collection of $n$ multiplicity sequences (that will be the multiplicity branches of a multiplicity tree). Since any multiplicity sequence is a sequence of integers that stabilizes to 1, we can describe them by the vectors
$$ M(i)=[m_{i,1},\ldots, m_{i,k_i}],$$
with the convention that $m_{i,j}=1$ for all $j \geq k_i-1$ and $m_{i,k_i-2}\neq 1$; it will be clear later why do not truncate the sequence to the last non-one entry.
\\
If $M=\max(k_1,\ldots,k_n)$ we write for all $i=1,\ldots,n$
$$ M(i)=[m_{i,1},\ldots, m_{i,M}],$$
in order to have vectors of the same length.
Each $M(i)$ represents a multiplicity sequence of an Arf numerical semigroup, so it must satisfy the following property:$$ \forall j\geq1 \textrm{ there exists }  s_{i,j} \in \mathbb{N}, \textrm{ such that }  s_{i,j}\geq j+1 \textrm{ and } m_{i,j}=\sum_{k=j+1}^{ s_{i,j}} m_{i,k}.$$
Denote by $\tau(E)$ the set of all multiplicity trees having the $n$ branches in $E$ and by $\sigma(E)$ the set of the corresponding Arf semigroups. We want to find an  unambiguous way to describe  distinct trees of $\tau(E)$.

 We define, for all $i=1,\ldots,n$, the following vectors
$$ S(i)=[ s_{i,1},\ldots s_{i,M}].$$
Because we have $m_ {i,j}=1$ for all $j\geq M-1$, it follows that $ s_{i,j}=j+1$ for all $j\geq M-1$.
\begin{example1} \label{ex1}
Let $M(1)$ be the following multiplicity sequence:
$$M(1)=[14,7,5,1,1]. $$ 
Then $S(1)$ is:
$$S(1)=[5,5,8,5,6]. $$  
Notice that, with this notation, from the vectors $S(i)$  we can easily reconstruct the sequences $M(i)$. It suffices to set  $m_{i,M}=1$  and then to compute the values of $m_{i,j}$ using the information contained in the integers $s_{i,j}$.   \end{example1} 

We will use the vectors  $S(i)$ to determine the level, in a tree of   $\tau(E)$, where two  branches have to split up.
\\
For each pair of integers $i,j$ such that $i<j$ and $i,j=1,\ldots,n$ we consider the  set $D(i,j)=\left\{ k : s_{i,k} \neq s_{j,k} \right\}$.
If $D(i,j) \neq \emptyset $  we consider the integer $$k_E(i,j)=\min\left\{ \min( s_{i,k}, s_{j,k}), k \in D(i,j) \right\}, $$
while if $D(i,j)=\emptyset$, and then the $i$-th and $j$-th branches have the same multiplcity sequence, we set $k_E(i,j)=+\infty$. We have the following proposition.

\begin{prop} \label{prop}
Consider a collection of multiplicity sequences $E$ and let $T \in \tau(E)$. Then $k_E(i,j)+1$ is the lowest level where the $i$-th and the $j$-th branches are prevented from being glued in $T$ (if $k_E(i,j)$ is infinite there are no limitations on the level where the branches have to split up).\end{prop} \noindent \textbf{Proof}
The case $k_E(i,j) = +\infty$ is trivial, because we have the same sequence along two consecutive branches and therefore no discrepancies that force the two branches to split up at a certain level.
Thus suppose $k_E(i,j) \neq +\infty$ and,  by contradiction, that the $i$-th and the $j$-th branches are glued at level $k_E(i,j)+1$.
 From the definition of $k_E(i,j)$, there exists $\overline{k} \in D(i,j)$ such that $k_E(i,j)=\min( s_{i,\overline{k}},s_{j,\overline{k}})$.
 Without loss of generality suppose that $\min( s_{i,\overline{k}},s_{j,\overline{k}})=s_{i,\overline{k}} \neq s_{j,\overline{k}}$ (since $\overline{k} \in D(i,j) $).

So in the tree we have the following nodes,
$$ (\ldots,m_{i,\overline{k}},\ldots,m_{j,\overline{k}},\ldots), \ldots,  (\ldots,m_{i,k_E(i,j)},\ldots,m_{j,k_E(i,j)},\ldots),$$$$, (\ldots,m_{i,k_E(i,j)+1},\ldots,m_{j,k_E(i,j)+1},\ldots).$$
We have that  $k_E(i,j)=s_{i,\overline{k}} $ so 
$$ m_{i,\overline{k}}=\sum_{t=\overline{k}+1}^{k_E(i,j)} {m_{i,t}},$$
while $k_E(i,j)+1= s_{i,\overline{k}}+1 \leq  s_{j,\overline{k}} $ so 
$$ m_{j,\overline{k}}\geq \sum_{t=\overline{k}+1}^{k_E(i,j)+1} {m_{j,t}}.$$
These facts easily imply that the first node cannot be expressed as a sum of the nodes of a subtree rooted in it, so we have a contradiction. 
Two branches are forced to split up only when we have this kind of problem, so the minimality of $k_E(i,j)$ 
guarantees that  they can be glued at level $k_E(i,j)$ (and obviously at lower levels). \qed
\begin{example2} \label{ex2}
Suppose that we have
$$M(1)=[14,7,5,1,1] \textrm{ and } M(2)=[7,3,1,1,1]. $$ 
So we have the vectors $S(1)$ and $S(2)$:
$$S(1)=[5,5,8,5,6] \textrm{ and } S(2)=[6,5,4,5,6]. $$ 
We have $D(1,2)=\left\{ 1,3 \right\}$, then $k(1,2)=\min \left\{ \min(5,6), \min(4,8) \right\}=\min\left\{5,4\right\}=4$.
Then the branches have to be separated at the fifth level.
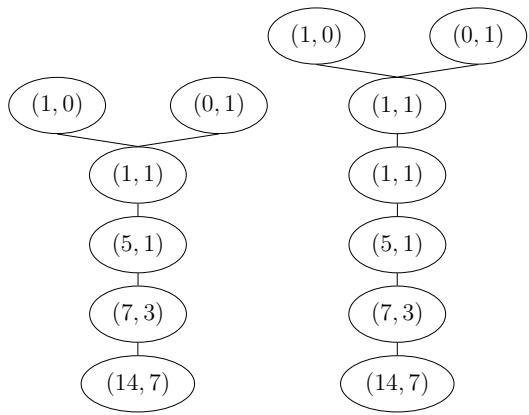
\begin{center}
\begin{tikzpicture}[grow'=up,sibling distance=32pt,scale=.75]
\tikzset{level distance=35pt,every tree node/.style={draw,ellipse}}
\Tree [.$(14,7)$ [.$(7,3)$ [.$(5,1)$ [.$(1,1)$ [.$(1,0)$ ] [.$(0,1)$ ] ] ] ] ]

\end{tikzpicture}\begin{tikzpicture}[grow'=up,sibling distance=32pt,scale=.75]
\tikzset{level distance=35pt,every tree node/.style={draw,ellipse}}\Tree [.$(14,7)$ [.$(7,3)$ [.$(5,1)$ [.$(1,1)$ [.$(1,1)$  [.$(1,0)$ ] [.$(0,1)$ ] ] ] ] ] ] \end{tikzpicture}
\end{center}

Notice that the first tree in the previous picture fulfills the properties of the multiplcity trees of an Arf semigroup. The second one cannot be the multiplicity tree of an Arf semigroup because the third node $(5,1)$ cannot be expressed as a sum of nodes in a subtree rooted in it.
\end{example2} 

Now we prove a general lemma that will be useful in the following.

\begin{Lemma1} \label{Lemma1}
Consider $v_1,v_2$ and $v_3$ in $\mathbb{N}^n$. If $i,j \in \left\{ 1,2,3 \right\}$  with $i \neq j$ we define:
\begin{itemize}
\item $\textrm{MIN}(v_i,v_j)=+\infty$ if $v_i=v_j$;
\item $\textrm{MIN}(v_i,v_j)=\min\left\{ \min(v_i[k],v_j[k]), k\in\left\{1,\ldots,n\right\}: v_i[k]\neq v_j[k]\right\}.$ 
\end{itemize}
Then there exists a permutation $\delta \in S^3$ such that $$ \textrm{MIN}(v_{\delta(1)},v_{\delta(2)})=\textrm{MIN}(v_{\delta(2)},v_{\delta(3)})\leq\textrm{MIN}(v_{\delta(1)},v_{\delta(3)}).$$ 
\end{Lemma1}
 \noindent \textbf{Proof} Suppose by contradiction that the thesis is not true. Then, renaming the indices  if necessary, we have
$$ \textrm{MIN}(v_1,v_2)<\textrm{MIN}(v_1,v_3)\leq\textrm{MIN}(v_2,v_3).$$
From the definition of $\textrm{MIN}(v_1,v_2)=l_{1,2}$ it follows that there exists a $k \in \left\{1,\ldots,n \right\}$ such that $v_1[k] \neq v_2[k]$ and $ \min(v_1[k],v_2[k])=l_{1,2}$. We have two cases:
\begin{itemize}
\item If $v_1[k]=l_{1,2}$ $\Rightarrow v_2[k] > l_{1,2}$. Then we must have $v_3[k]=l_{1,2}$, in fact otherwise we would have  $\textrm{MIN}(v_1,v_3)\leq l_{1,2}<\textrm{MIN}(v_1,v_3)$.
But then$$  l_{1,2}<\textrm{MIN}(v_2,v_3) \leq \min(v_2[k],v_3[k])=l_{1,2},$$
and we have a contradiction.
\item If $v_2[k]=l_{1,2}$ $\Rightarrow v_1[k] > l_{1,2}$. Then we must have $v_3[k]=l_{1,2}$, in fact otherwise we would have  $\textrm{MIN}(v_2,v_3)\leq l_{1,2}<\textrm{MIN}(v_2,v_3)$.
But then$$  l_{1,2}<\textrm{MIN}(v_1,v_3) \leq \min(v_1[k],v_3[k])=l_{1,2},$$
and we have a contradiction.\qed
\end{itemize} 
\begin{rem} \label{rem}
If we have three multiplicity sequences $M(1)$, $M(2)$ and $M(3)$ then, if $E=\left\{M(1),M(2),M(3) \right\}$  then there exist a permutation $\delta \in S^3$  such that
 $$ k_E(\delta(1),\delta(2))= k_E(\delta(2),\delta(3)) \leq  k_E(\delta(1),\delta(3)).$$
In fact the integers $k_E(i,j)$ are of the same type of the integers $\textrm{MIN}(v_i,v_j)$ of the previous lemma with $v_i=S(i)$.
\end{rem}

We give now a way to describe a tree of $\tau(E)$.
 If $T \in \tau(E)$, it can be represented by an upper triangular matrix $ n \times n$  $$M(T)_{E}=\left(  \begin{matrix} 0 & p_{1,2} & p_{1,3} & \ldots & p_{1,n} \\ 
0 & 0 & p_{2,3} & \ldots & p_{2,n} \\ \ldots & \ldots & \ldots & \ldots & \ldots \\ 0 & 0 & 0 & \ldots& p_{n-1,n} \\ 0 & 0 & 0 & \ldots & 0 \\  \end{matrix}\right),$$ where $p_{i,j}$ is the highest level such that the  $i$-th and the $j$-th branches are glued in $T$.
\begin{remark} \label{remark}
If $M(T)_{E}$ is the matrix of a $T \in \tau(E)$, it is clear that everytime we consider three indices $i<j<k$ we must have:
$$ p_{i,j} \geq \min(p_{i,k},p_{j,k}), p_{j,k} \geq \min(p_{i,j},p_{i,k}) \textrm{ and } p_{i,k} \geq \min(p_{i,j},p_{j,k}),$$
when we are using the obvious fact that the relation of being glued has the transitive property.
From the previous inequalities it follows that the set $\left\{ p_{i,j}, p_{j,k}, p_{i,k}\right\}=\left\{ x,x,y \right\}$, with $x\leq y$ (independently of the order).
\end{remark}
 From Proposition \ref{prop} we have that $p_{i,j} \in \left\{ 1,\ldots, k_E(i,j) \right\}$ for all $i,j=1,\ldots,n$ with $i<j$. 
In the following, with an abuse of notation, we will identify a tree with its representation.

We call a tree $T$ of $\tau(E)$ untwisted if two nonconsecutive branches are glued at level $l$ if and only if  all the consecutive branches  between them are glued  at a level greater or equal to $l$. We will call twisted a tree that it is not untwisted.

From the definition it follows that the matrix of an untwisted tree $T \in \tau(E)$ is such that: $$ p_{i,j}=\min\left\{ p_{i,i+1},p_{i+1,i+2},\ldots,p_{j-1,j}\right\} \textrm{ for all } i<j. $$
So an untwisted tree can be completely described by the second diagonal of its matrix. Thus in the following we will  indicate an untwisted tree by a vector $T_E=(d_1,\ldots,d_{n-1})$ where $d_i=p_{i,i+1}$.
It is easy to see that a twisted tree can be converted to an untwisted one by accordingly permuting its branches.
So in the following we can focus, when it is possible, only on the properties of the untwisted trees, that are easier to study, obtaining the twisted one by permutation. 

\begin{example23}
Let us consider the following tree of $\tau(E)$ with $$E=\left\{ M(1)=[5,4,1,1],M(2)=[2,2,1,1],M(3)=[6,4,1,1]\right\}
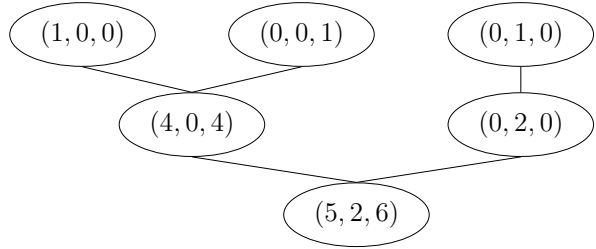
$$
\begin{center}
\begin{tikzpicture}[grow'=up,sibling distance=32pt,scale=.85]
\tikzset{level distance=40pt,every tree node/.style={draw,ellipse}}
\Tree [ .$(5,2,6)$ [ .$(4,0,4)$ [ .$(1,0,0)$ ] [ .$(0,0,1)$ ] ] [ .$(0,2,0)$ [ .$(0,1,0)$ ] ] ]

\end{tikzpicture}
\end{center}

This tree is twisted because the first and the third branches are glued at level two while the first and the second are not.

If we consider the permutation $(2,3)$ on the branches we obtain the tree

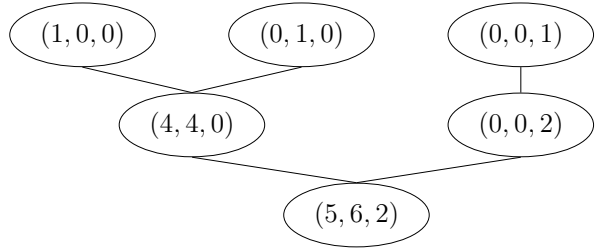
\begin{center}
\begin{tikzpicture}[grow'=up,sibling distance=32pt,scale=.85]
\tikzset{level distance=40pt,every tree node/.style={draw,ellipse}}
\Tree [ .$(5,6,2)$ [ .$(4,4,0)$ [ .$(1,0,0)$ ] [ .$(0,1,0)$ ] ] [ .$(0,0,2)$ [ .$(0,0,1)$ ] ] ]

\end{tikzpicture}
\end{center}
that is untwisted, even if it belongs to a different set $\tau(E')$ where  $$E'=\left\{ M(1)=[5,4,1,1],M(2)=[6,4,1,1],M(3)=[2,2,1,1]\right\},$$ and can be represented by the vector $T_E'=(2,1)$.
\end{example23}
Denote by $S(T)$ the semigroup determined by the tree $T$. In  \cite[Lemma 5.1]{BDF2} it is shown that if $T^1$ and $T^2$ are untwisted trees of $\tau(E)$, then $S(T^1) \subseteq S(T^2) $ if and only if    $T^2_{E}$ is coordinatewise less than or equal to   $T^1_{E}$. The previous result can be easily extended to the twisted trees. Then, in the general case we have that  $S(T^1) \subseteq S(T^2) $, where $S(T^1)$ and $S(T^2)$ belong to $\sigma(E)$, if and only if each entry of  $M(T^2)_{E}$ is  less than or equal to the corresponding entry of   $M(T^1)_{E}$.
If $k_E(i,j) \neq +\infty$ for all $i<j$, we can  consider  $T^{\textrm{MIN}}$ such that $$M(T^{\textrm{MIN}})_{E}=\left(  \begin{matrix} 0 & k_E(1,2) & k_E(1,3) & \ldots & k_E(1,n) \\ 
0 & 0 & k_E(2,3) & \ldots & k_E(2,n) \\ \ldots & \ldots & \ldots & \ldots & \ldots \\ 0 & 0 & 0 & \ldots& k_E(n-1,n) \\ 0 & 0 & 0 & \ldots & 0 \\  \end{matrix}\right),$$ 
that is well defined for Remark \ref{rem}. Then  $S(T^{\textrm{MIN}})$ is the smallest Arf semigroup belonging to $\sigma(E)$.
\begin{remark1}
If in the collection $E$ there are two  branches with the same multiplicity sequence then $ |\sigma(E)|=+ \infty$.
\end{remark1}
\begin{example3}
\label{ex3}
We can  count the number of  untwisted trees of $\tau(E)$ by using their representation. If we call $\tau^*(E)$ the set of all the untwisted trees of $\tau(E)$, these trees are completely determined by the elements in the second diagonal of their matrix, that are bounded by $k_E(j,j+1)$.
Hence the number of untwisted trees is:
$$ |\tau^*(E)|= \prod_{j=1}^{n-1}{k_E(j,j+1)}.$$
Suppose that $E= \left\{ M(1),M(2),M(3) \right\}$, where
$$ M(1)=[5,4,1,1]\textrm{, }  M(2)=[6,4,1,1]\textrm{, } M(3)=[2,2,1,1]. $$
We have:
$$ S(1)=[3,6,4,5]\textrm{, } S(2)=[4,6,4,5]\textrm{, } S(3)=[2,4,4,5]. $$
Then $D(1,2)=\left\{1\right\}, D(2,3)=\left\{1,2\right\}$ and $k_E(1,2)=\min(3,4)=3$ and \newline $k_E(2,3)=\min\left\{ \min(2,4),\min(4,6)\right\}=2$. 
There are $k_E(1,2) \cdot k_E(2,3)=6$ trees in $\tau^*(E)$.
 They are:
\vskip 0.3in

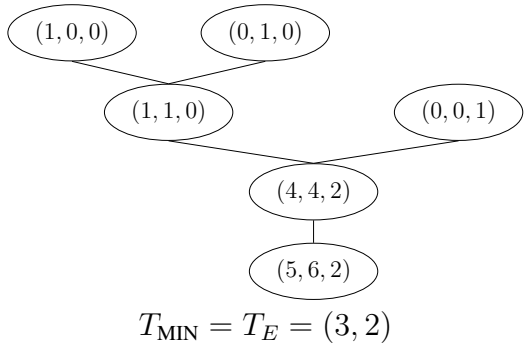
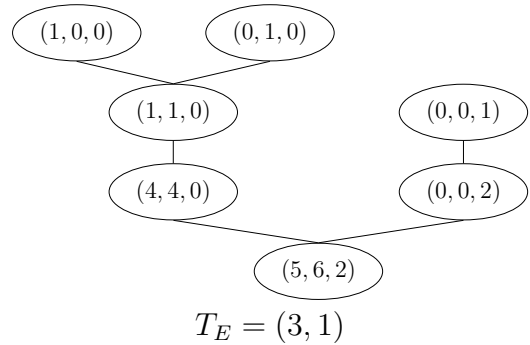
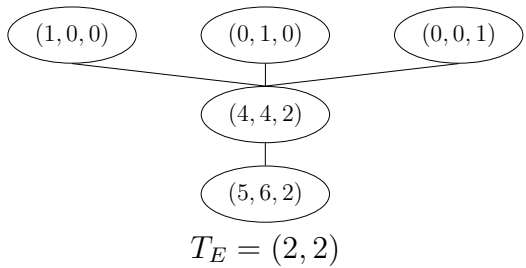
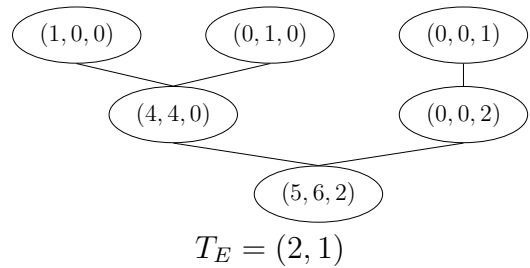
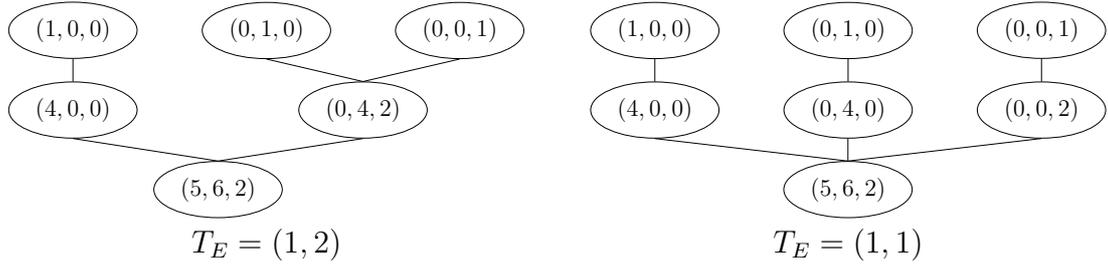
\begin{tikzpicture}[grow'=up,sibling distance=32pt,scale=.75]
\tikzset{level distance=40pt,every tree node/.style={draw,ellipse}} \Tree [ .$(5,6,2)$ [ .$(4,4,2)$ [ .$(1,1,0)$ [ .$(1,0,0)$ ] [ .$(0,1,0)$ ] ] [ .$(0,0,1)$ ]] ] \node[below]at(current bounding box.south){$T_{\textrm{MIN}}=T_E=(3,2)$}; \end{tikzpicture}
\hskip 0.3in
\begin{tikzpicture}[grow'=up,sibling distance=32pt,scale=.75]
\tikzset{level distance=40pt,every tree node/.style={draw,ellipse}} \Tree [ .$(5,6,2)$ [ .$(4,4,0)$ [ .$(1,1,0)$ [ .$(1,0,0)$ ] [ .$(0,1,0)$ ] ] ] [ .$(0,0,2)$ [ .$(0,0,1)$ ] ] ]\node[below]at(current bounding box.south){$T_E=(3,1)$}; \end{tikzpicture}
\vskip 0.3in
\begin{tikzpicture}[grow'=up,sibling distance=32pt,scale=.75]
\tikzset{level distance=40pt,every tree node/.style={draw,ellipse}} \Tree [ .$(5,6,2)$ [ .$(4,4,2)$ [ .$(1,0,0)$ ] [ .$(0,1,0)$ ] [ .$(0,0,1)$ ] ] ] \node[below]at(current bounding box.south){$T_E=(2,2)$}; \end{tikzpicture}
\hskip 0.3in
\begin{tikzpicture}[grow'=up,sibling distance=32pt,scale=.75]
\tikzset{level distance=40pt,every tree node/.style={draw,ellipse}} \Tree [ .$(5,6,2)$ [ .$(4,4,0)$ [ .$(1,0,0)$ ] [ .$(0,1,0)$ ] ] [ .$(0,0,2)$ [ .$(0,0,1)$ ] ] ]\node[below]at(current bounding box.south){$T_E=(2,1)$}; \end{tikzpicture}
\vskip 0.3in
\begin{tikzpicture}[grow'=up,sibling distance=32pt,scale=.75]
\tikzset{level distance=40pt,every tree node/.style={draw,ellipse}} \Tree [ .$(5,6,2)$ [ .$(4,0,0)$ [ .$(1,0,0)$ ] ] [ .$(0,4,2)$ [ .$(0,1,0)$ ] [ .$(0,0,1)$ ] ] ]\node[below]at(current bounding box.south){$T_E=(1,2)$}; \end{tikzpicture}
\hskip 0.3in
\begin{tikzpicture}[grow'=up,sibling distance=32pt,scale=.75]
\tikzset{level distance=40pt,every tree node/.style={draw,ellipse}} \Tree [ .$(5,6,2)$  [ .$(4,0,0)$ [.$(1,0,0)$ ] ] [ .$(0,4,0)$  [.$(0,1,0)$ ]] [ .$(0,0,2)$ [.$(0,0,1)$ ] ] ]\node[below]at(current bounding box.south){$T_E=(1,1)$}; \end{tikzpicture}

\begin{rem1}
Because we are able to determine completely $\tau^*(E)$ for each $E$ collection of multiplicity sequences we have a way to determine $\tau(E)$. If $\delta \in S_n$ is a permutation of the symmetric group $S_n$ we can consider $\delta^{-1}(\tau^*(\delta(E))) \subseteq \tau(E).$ It is trivial to see that
$$ \bigcup_{\delta \in S_n}{\delta^{-1}(\tau^*(\delta(E)))}=\tau(E).$$
\end{rem1}
If we apply this strategy to find $\tau(E)$ with the $E$ of the previous example we find that in $\tau(E)$ there is only one twisted tree $T$ with $$ M(T)_E=\left( \begin{matrix} 0 &1 &2 \\ 0& 0 &1 \\ 0 &0 &0 \\ \end{matrix}\right). $$ \vskip 0.3in
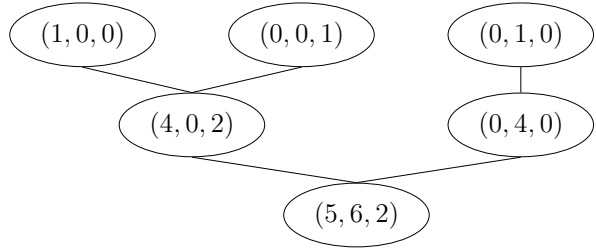
\begin{center}\begin{tikzpicture}[grow'=up,sibling distance=32pt,scale=.85]
\tikzset{level distance=40pt,every tree node/.style={draw,ellipse}} \Tree [ .$(5,6,2)$ [ .$(4,0,2)$ [ .$(1,0,0)$ ] [ .$(0,0,1)$ ] ] [ .$(0,4,0)$ [ .$(0,1,0)$ ] ] ] \end{tikzpicture} \end{center}

\end{example3}

\section{ When a set of vectors determines an Arf semigroup} \label{section3}
In this section we want to understand when a set $G \subseteq \mathbb{N}^n$ determines uniquely an Arf semigroup of  $\mathbb{N}^n$.
 First of all we need to fix some notations.
 
Given $ G\subseteq \mathbb{N}^n$ we denote by $S(G)$ the following set $$S(G)=\left\{ S: S \subseteq \mathbb{N}^n \textrm{ is an Arf semigroup and } G \subseteq S \right\}.$$
If the set $ S(G)$ has a minimum (with the partial order given by the inclusion), we will denote such a minimum by $\textrm{Arf}(G)$. 
Hence we have to understand when $\textrm{Arf}(G)$ is well defined and, in this case, how to find it.

 If $i \in \left\{ 1,\ldots,n \right\}$, and $S \in S(G)$ we denote by $S_i$ the projection on the $i$-th coordinate. We know that $S_i$ is an Arf numerical semigroup and it contains the set $G[i]=\left\{ g[i]: g \in G \right\}$ where with $g[i]$ we indicate the $i$-th coordinate of the vector $g$.
We recall also that, if we have a set of integers $I$ such that $\gcd(I)=1$, it is possible to compute the smallest Arf semigroup containing $I$, that is the Arf closure of the numerical semigroup generated by the elements of $I$.
This computation can be made by using  the modified Jacobian algorithm of Du Val (cf. \cite{Duval}). 

We have the following theorem:

\begin{theo}
\label{T45}
Suppose that we have $G \subseteq \mathbb{N}^n$.  Then $\textrm{Arf}(G)$ is well defined if and only if the following conditions hold:
\begin{itemize}
\item $\gcd\left\{ g[i], g \in G \right\}=1$ for $i=1,\ldots,n;$
\item For all $i,j \in \left\{1,\ldots,n \right\}$ such that $i < j$ there exists  $g \in G$ such that $g[i] \neq g[j]$. 
\end{itemize}
\end{theo}

\noindent \textbf{ Proof} ($\Rightarrow $)
Suppose that   $\textrm{Arf}(G)$ is well defined and suppose by contradiction that the two  conditions of the theorem are not simultaneuously fulfilled.

We have two cases.

\begin{itemize}
\item \textbf{Case 1}: The first condition is not fulfilled.

Then there exists an $i $ such that  $\gcd(G[i])=d \neq 1$. When we apply the Jacobian algorithm to the elements of $G[i]$ we will produce a sequence of the following type:
$$ [m_{i,1}, \ldots, m_{i,k},\ldots]$$
where there exists a $k$  such that $m_{i,j}=d$ for all $j \geq k$ (it happens because the Jacobian algorithm performs an Euclidean algorithm on $G[i]$). Denote by $\overline{k}$ the minimum $k$ such that the Arf semigroup associated to the sequence $$ [m_{i,1}, \ldots, m_{i,\overline{k}}=d,1,1],$$ contains $G[i]$ (such minimum exists for the properties of the algorithm of Du Val). Then for all $z \geq \overline{k}$ we can consider the multiplicity sequence
$$ M(z)=[m_{i,1}, \ldots, m_{i,\overline{k}}=d,\ldots, m_{i,z}=d,1,1]$$
and if $S(z)$ is the Arf numerical semigroup associated to $M(z)$ then $ G[i] \subseteq S(z)$.
 Now it is trivial to show that $S(z_1) \subseteq S(z_2)$ if $z_1 \geq z_2$. 
Then we have an infinite decreasing chain of Arf semigroup containing the set $G[i]$. This means that the projection on the $i$-th branch can be smaller and smaller, therefore we cannot find a minimum in the set $S(G)$.

Thus we have found a contradiction in this case.

An example illustrating \textbf{Case 1} is the following.

If we consider $G=\left\{ [2,3],[4,4] \right\} $, we have no information on the multiplicity sequence along the first branch and so we can obtain the following infinite decreasing chain of Arf semigroups containing $G$:

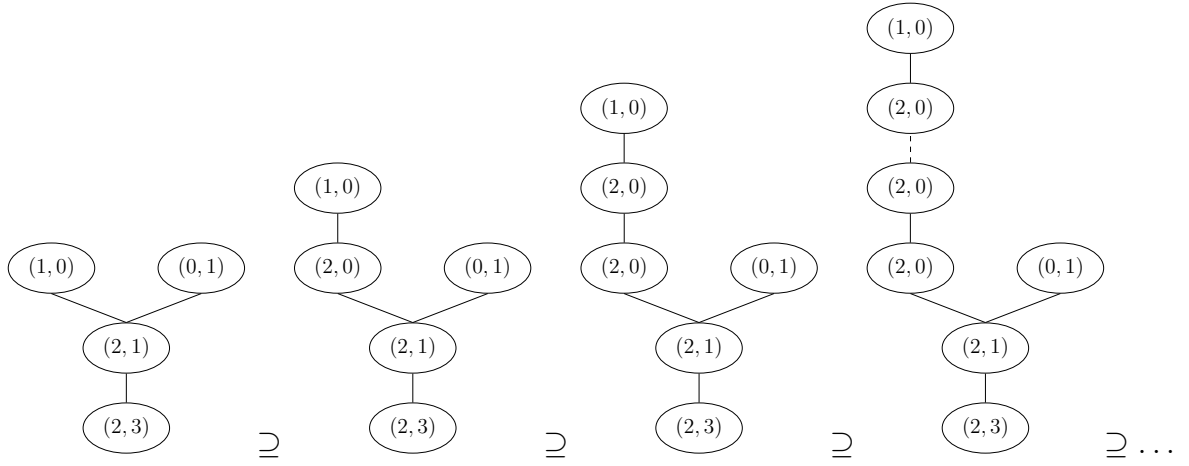
\begin{tikzpicture}[grow'=up,sibling distance=35pt,scale=.67]
\tikzset{level distance=45pt,every tree node/.style={draw,ellipse}} \Tree [ .$(2,3)$ [ .$(2,1)$ [ .$(1,0)$ ] [ .$(0,1)$ ] ] ]\end{tikzpicture}
\begin{centering} $\supseteq$ \end{centering}  \begin{tikzpicture}[grow'=up,sibling distance=35pt,scale=.67]
\tikzset{level distance=45pt,every tree node/.style={draw,ellipse}} \Tree [ .$(2,3)$ [ .$(2,1)$ [ .$(2,0)$ [ .$(1,0)$ ] ] [ .$(0,1)$ ] ] ]\end{tikzpicture} 
\begin{centering} $\supseteq$ \end{centering}  \begin{tikzpicture}[grow'=up,sibling distance=35pt,scale=.67]
\tikzset{level distance=45pt,every tree node/.style={draw,ellipse}} \Tree [ .$(2,3)$ [ .$(2,1)$ [ .$(2,0)$ [ .$(2,0)$ [ .$(1,0)$ ]] ] [ .$(0,1)$ ] ] ]\end{tikzpicture} \begin{centering} $\supseteq$ \end{centering}  \begin{tikzpicture}[grow'=up,sibling distance=35pt,scale=.67]
\tikzset{level distance=45pt,every tree node/.style={draw,ellipse}} \Tree [ .$(2,3)$ [ .$(2,1)$ [ .$(2,0)$ [ .$(2,0)$ \edge[dashed,level distance=16pt];[ .$(2,0)$ [  .$(1,0)$ ] ]] ] [ .$(0,1)$ ] ] ]\end{tikzpicture}$\supseteq \nobreak\dots$

\item \textbf{Case 2}: The first condition is  fulfilled.

So in this case the second condition is not fulfilled.
The fact that $\gcd\left\{ g[i], g \in G \right\}=1$ for $i=1,\ldots,n$ implies that we can compute the smallest Arf numerical semigroup $S(i)$ containing $G[i]$ for all $i=1, \ldots,n$.

Therefore if we denote by $M_i$ the multiplicity sequence of $S(i)$ we clearly have that $\textrm{Arf}(G) \in \sigma(E), $ where $E=\left\{ M_i, i=1,\ldots,n \right\}$. Suppose that it is defined by the matrix  $$M(T)_{E}=\left(  \begin{matrix} 0 & p_{1,2} & p_{1,3} & \ldots & p_{1,n} \\ 
0 & 0 & p_{2,3} & \ldots & p_{2,n} \\ \ldots & \ldots & \ldots & \ldots & \ldots \\ 0 & 0 & 0 & \ldots& p_{n-1,n} \\ 0 & 0 & 0 & \ldots & 0 \\  \end{matrix}\right).$$

Now if we consider an element $ h \in G[i]$ we have that $h \in S(i)$ and therefore there exists an index $\textrm{pos}_E(i,h)$ such that
$$ h= \sum_{k=1}^{\textrm{pos}_E(i,h)} {m_{i,k}}.$$
 If $g\in G$ we can define $\textrm{pos}_E(g)=[\textrm{pos}_E(1,g[1]),\ldots, \textrm{pos}_E(n,g[n])]$.

Notice that, if we consider $i,j \in \left\{ 1,\ldots,n \right\}$, with $i<j$ and  $g \in G$ such that $\textrm{pos}_E(i,g[i]) \neq \textrm{pos}_E(j,g[j])$, we can easily deduce that in a multiplicity tree of an  Arf semigroup of $\sigma(E)$ containing $G$ the $i$-th and $j$-th branches cannot be glued at a level greater than $\min(\textrm{pos}_E(i,g[i]), \textrm{pos}_{E}(j,g[j]))$ . 

Then $p_{i,j}$ is at most $\min(\textrm{pos}_E(i,g[i]), \textrm{pos}_{E}(j,g[j]))$, and we also have to recall that $p_{i,j}$ is at most $k_E(i,j)$.

So denote by $$U_E(G)=\left\{(i,j) \in \left\{ 1,\ldots, n \right\}^2 : i<j; \textrm{pos}_E(i,g[i])=\textrm{pos}_{E}(j,g[j]) \textrm{ for all } g \in G \right\}.$$
For each $ (i,j) \notin U_E(G)$ we define
$$ \textrm{MIN}_E(i,j,G)=\min\left(k_E(i,j),\min\left\{\min(\textrm{pos}_E(i,g[i]), \textrm{pos}_E(j,g[j])):  g \in G, \right. \right. $$ $$ \left. \left.  \textrm{pos}_E(i,g[i]) \neq \textrm{pos}_{E}(j,g[j]) \right\} \right) .$$
Notice that we need $(i,j) \notin U_E(G)$ to have the previous integers well defined.

So from the previous remark it follows that an Arf semigroup $S(T^1)$  of $\sigma(E)$ containing $G$ with $$M(T^1)_{E}=\left(  \begin{matrix} 0 & a_{1,2} & a_{1,3} & \ldots & a_{1,n} \\ 
0 & 0 & a_{2,3} & \ldots & a_{2,n} \\ \ldots & \ldots & \ldots & \ldots & \ldots \\ 0 & 0 & 0 & \ldots& a_{n-1,n} \\ 0 & 0 & 0 & \ldots & 0 \\  \end{matrix}\right) $$ is such that $a_{i,j}$ is at most $k_E(i,j)$ for $(i,j) \in U_E(G)$ and $a_{i,j}$ is at most $\textrm{MIN}_E(i,j,G)$ for $(i,j) \notin U_E(G)$.
So for the Arf closure we want to choose the biggest possible values, therefore we have:
$$ p_{i,j}=k_E(i,j)  \textrm{ for } (i,j) \in U_E(G) \textrm{  and  } p_{i,j}=\textrm{MIN}_E(i,j,G) \textrm{ for } (i,j) \notin U_E(G).$$
We need to prove that this integers are compatible with the transitive property of a matrix of an Arf semigroup tree. Therefore we consider a triad of integers $i<j<k$ and we want to show that $p_{i,j},p_{j,k}$ and $p_{i,k}$ are in a $ \left\{x,x,y\right\}$ configuration. We have the following cases:
\begin{enumerate}
\item $(i,j),(j,k),(k,i) \in U_E(G)$. Then $p_{i,j}=k_E(i,j),$$p_{i,k}=k_E(i,k)$ and $p_{j,k}=k_E(j,k)$ and for the Remark \ref{rem} they satisfy our condition;
\item $(i,j),(j,k),(k,i) \notin U_E(G)$. We consider the vectors $$v_l=[\textrm{pos}_E(l,g_1[l]),\ldots,\textrm{pos}_E(l,g_m[l])],$$
where $l\in \left\{i,j,k\right\}$ and $G=\left\{g_1,\ldots,g_m\right\}$. 
Then, using the notations of  Lemma \ref{Lemma1}, we have that $$ p_{i,j}=\min(k_E(i,j),\textrm{MIN}(v_i,v_j)),p_{i,k}=\min(k_E(i,k),\textrm{MIN}(v_i,v_k)) \textrm{ and } $$$$p_{j,k}=\min(k_E(j,k),\textrm{MIN}(v_j,v_k)).$$
Then suppose by contradiction that they are not compatible. Without loss of generality we can assume that $$ p_{i,j}<p_{i,k} \leq p_{j,k}. $$
We have two cases 
\begin{itemize}
\item $p_{i,j}=k_E(i,j).$  Then we would have
$$ k_E(i,j)=p_{i,j}<p_{j,k} \leq k_E(j,k) \textrm{ and } k_E(i,j)=p_{i,j}<p_{i,k} \leq k_E(i,k), $$
and this is absurd for the Remark \ref{rem};
\item$p_{i,j}=\textrm{MIN}(v_i,v_j).$  Then we would have
$$ \textrm{MIN}(v_i,v_j)=p_{i,j}<p_{j,k} \leq\textrm{MIN}(v_j,v_k)\textrm{ and } \textrm{MIN}(v_i,v_j)=p_{i,j}<p_{i,k} \leq\textrm{MIN}(v_i,v_k), $$
and this is absurd  against  Lemma \ref{Lemma1} applied to the vectors $v_i,v_j$ and $v_k$.
\end{itemize}
\item $(i,j) \in U_E(G)$ and $(j,k),(k,i) \notin U_E(G)$ (and the similar configurations).
In this case we have that $v_i=v_j$. Then
$$ p_{i,j}=k_E(i,j), p_{i,k}=\min(k_E(i,k),x), \textrm{ and }   p_{j,k}=\min(k_E(j,k),x),$$
where $x=\textrm{MIN}(v_i,v_k)= \textrm{MIN}(v_j,v_k)$.
We have two cases:
\begin{itemize}
\item $k_E(i,j)=k_E(j,k) \leq k_E(i,k) $ (or equivalently $k_E(i,j)=k_E(i,k) \leq k_E(j,k) $).
If $x<k_E(j,k)\leq k_E(i,k)$ then we have $ p_{j,k}=p_{i,k}=x<k_E(i,j)$ and it is fine.
If $x \geq k_E(j,k)$ then $p_{j,k}=k_E(j,k)=p_{i,j} \leq p_{i,k}$ that is compatible too.
\item $k_E(i,k)=k(j,k) < k_E(i,j)$. In this case we have $p_{i,k}=p_{j,k}<k_E(i,j)=p_{i,j}$ and it is fine.
\end{itemize}
\end{enumerate}

 So we actually have a well defined tree.

Anyway, because the second condition is not fulfilled, then there exists a pair $ (i,j)\in \left\{1,\ldots,n \right\}^2$ such that for all $g \in G$ we have $g[i] = g[j]$.  So $(i,j) \in U_E(G)$, and, since in this case the two sequences are the same, we  obtain $p_{i,j}=k_E(i,j)=+\infty$.

Thus we have found a contradiction because $\textrm{Arf}(G)$ is not well defined.

An example illustrating \textbf{Case 2} is the following.
If we consider $G=\left\{ [3,3,2],[2,2,3] \right\}$, we will have the same multiplicity sequences in the first two branches, with no clues about the splitting point so  we can obtain the following infinite decreasing chain in $S(G)$:
\vskip 0.3in
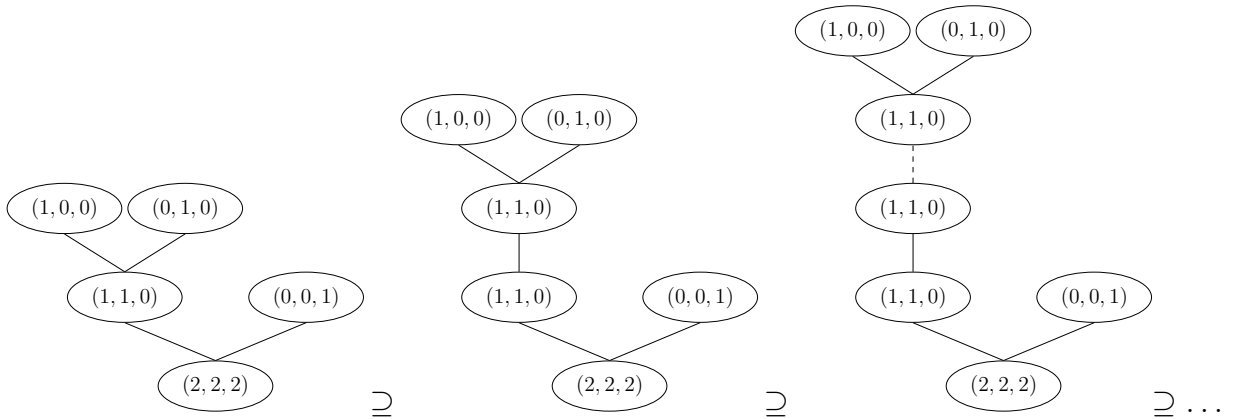
\begin{tikzpicture}[grow'=up,sibling distance=3pt,scale=.67]
\tikzset{level distance=50pt,every tree node/.style={draw,ellipse}} \Tree [ .$(2,2,2)$ [ .$(1,1,0)$ [ .$(1,0,0)$ ] [ .$(0,1,0)$ ] ] [ .$(0,0,1)$ ] ]\end{tikzpicture}
 $\supseteq$  \begin{tikzpicture}[grow'=up,sibling distance=3pt,scale=.67]
\tikzset{level distance=50pt,every tree node/.style={draw,ellipse}} \Tree [ .$(2,2,2)$ [ .$(1,1,0)$ [ .$(1,1,0)$ [ .$(1,0,0)$ ] [ .$(0,1,0)$ ] ] ] [ .$(0,0,1)$ ] ]\end{tikzpicture}
$\supseteq$    \begin{tikzpicture}[grow'=up,sibling distance=3pt,scale=.67]
\tikzset{level distance=50pt,every tree node/.style={draw,ellipse}} \Tree [ .$(2,2,2)$ [ .$(1,1,0)$ [ .$(1,1,0)$ \edge[dashed] ;[ .$(1,1,0)$ [ .$(1,0,0)$ ] [ .$(0,1,0)$ ] ] ] ] [ .$(0,0,1)$ ] ]\end{tikzpicture}\nobreak $\supseteq \nobreak\dots$

\end{itemize}

\noindent ($\Leftarrow$)
The previous proof gives us  a way to compute $\textrm{Arf}(G)$.
We have  to compute, using the modified Jacobian algorithm of Du Val, the  Arf closure of each $G[i]$, finding the collection $E$ (the first condition guarantees that it is possible to do that). Then we can find the matrix describing the semigroup using the set $U_E(G)$ and the integers $\textrm{MIN}_E(i,j,G)$ with the procedure present in the first part (we cannot have $p_{i,j}=+\infty$ for the second condition).  \qed

\begin{example6}\label{ex6}

 Suppose that we have $ G=\left\{\textrm{G}(1)=[5,6,5], \textrm{G}(2)=[9,11,4], \textrm{G}(3)=[9,10,2] \right\}, $
 that satisfies the conditions of the theorem.
Then we have to apply the modified Jacobian algorithm to the sets
 $$G[1]=\left\{5,9\right\}, G[2]=\left\{6,10,11 \right\} \textrm{ and } G[3]=\left\{  2,4,5\right\}.$$
We will find the following multiplicity sequences:
$$ M_1=[5,4,1,1],  M_2=[6,4,1,1] \textrm{ and } M_3=[2,2,1,1]. $$
We have $k_E(1,2)=3$, $k_E(2,3)=2$ and $k_E(1,3)=2$.

So we have $ \textrm{pos}_E(\textrm{G}(1))=[1,1,3]$,  $ \textrm{pos}_E(\textrm{G}(2))=[2,3,2]$ and  $ \textrm{pos}_E(\textrm{G}(3))=[2,2,1].$

In this case $U_E(G)=\emptyset$. 

We have $ \textrm{MIN}_E(1,2,G)=\min(2,k_E(1,2))=2$, $ \textrm{MIN}_E(2,3,G)=\min(1,k_E(2,3))=1$ and $ \textrm{MIN}_E(1,3,G)=\min(1,k_E(1,3))=1$.

So the Arf closure is described by the matrix $$ M(T)_E=\left( \begin{matrix} 0 &2 &1 \\ 0& 0 &1 \\ 0 &0 &0 \\ \end{matrix}\right). $$ with $$E=\left\{M_1=[5,4,1,1],  M_2=[6,4,1,1] \textrm{ and } M_3=[2,2,1,1] \right\}.$$
Notice that in this case we find that the Arf closure is an untwisted tree of $\tau(E)$ represented by the vector $T_E=(2,1)$.

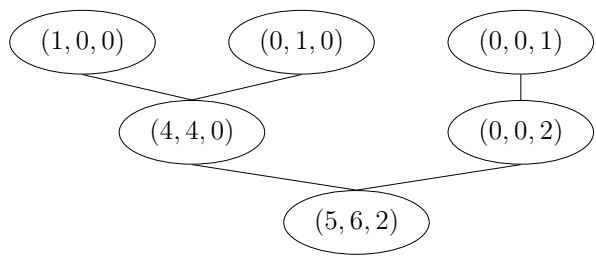
\begin{center}
\begin{tikzpicture}[grow'=up,sibling distance=32pt,scale=.85]
\tikzset{level distance=40pt,every tree node/.style={draw,ellipse}} \Tree [ .$(5,6,2)$ [ .$(4,4,0)$ [ .$(1,0,0)$ ] [ .$(0,1,0)$ ] ] [ .$(0,0,2)$ [ .$(0,0,1)$ ] ] ] \end{tikzpicture} \end{center}

\end{example6}
\begin{example67}
 Suppose that we have $ G=\left\{\textrm{G}(1)=[8,6,10], \textrm{G}(2)=[5,10,5], \textrm{G}(3)=[10,13,8] \right\}, $
 that satisfies the conditions of the theorem.
Then we have to apply the modified Jacobian algorithm to the sets
 $$G[1]=\left\{5,8,10\right\}, G[2]=\left\{6,10,13 \right\} \textrm{ and } G[3]=\left\{ 5,8, 10\right\}.$$
We will find the following multiplicity sequences:
$$ M_1=[5,3,2,1,1],  M_2=[6,4,2,1,1] \textrm{ and } M_3=[5,3,2,1,1]. $$
We have $k_E(1,2)=4$, $k_E(2,3)=4$ and $k_E(1,3)=+ \infty $.

So we have $ \textrm{pos}_E(\textrm{G}(1))=[2,1,3]$,  $ \textrm{pos}_E(\textrm{G}(2))=[1,2,1]$ and  $ \textrm{pos}_E(\textrm{G}(3))=[3,4,2].$

In this case $U_E(G)=\emptyset$. 

We have $ \textrm{MIN}_E(1,2,G)=\min(1,k_E(1,2))=1$, $ \textrm{MIN}_E(2,3,G)=\min(1,k_E(2,3))=1$ and $ \textrm{MIN}_E(1,3,G)=\min(2,k_E(1,3))=2$.

So the Arf closure is described by the matrix $$ M(T)_E=\left( \begin{matrix} 0 &1 &2\\ 0& 0 &1 \\ 0 &0 &0 \\ \end{matrix}\right). $$ with $$E=\left\{M_1=[5,3,2,1,1],  M_2=[6,4,2,1,1] \textrm{ and } M_3=[5,3,2,1,1] \right\}.$$
Notice that in this case we find that the Arf closure is a twisted tree.
\vskip 0.3in
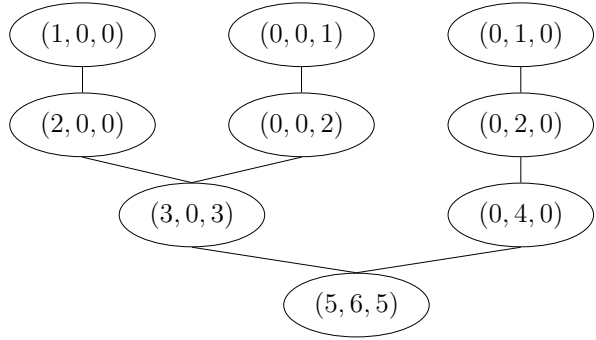
\begin{center}
\begin{tikzpicture}[grow'=up,sibling distance=32pt,scale=.85]
\tikzset{level distance=40pt,every tree node/.style={draw,ellipse}} \Tree [ .$(5,6,5)$ [ .$(3,0,3)$ [ .$(2,0,0)$ [ .$(1,0,0)$ ] ] [ .$(0,0,2)$ [ .$(0,0,1)$ ] ] ] [ .$(0,4,0)$ [ .$(0,2,0)$ [ .$(0,1,0)$ ]] ]] \end{tikzpicture} \end{center}

\end{example67}

\section{Arf closure of a good semigroup of $\mathbb{N}^n$} \label{section4}

Denote by $S$ a good semigroup of $\mathbb{N}^n$. In this section we describe a way to find the smallest Arf semigroup of $\mathbb{N}^n$ containing $S$, that is the Arf closure of $S$ (the existence of the Arf closure is proved in \cite{DPMT}). We denote this semigroup by $\textrm{Arf}(S)$.
If $S$ is a good semigroup of $\mathbb{N}^n$, we denote by $S_i$ the projection on the $i$-th coordinate. The properties of the good semigroups guarantee that $S_i$ is a numerical semigroup. Thus it is clear that an Arf semigroup $T$ containing $S$ is such that $\textrm{Arf}(S_i) \subseteq T_i$ for all $i=1,\ldots,n$, where $\textrm{Arf}(S_i)$ is the Arf closure of the numerical semigroup $S_i$  (we can compute it using the algorithm of Du Val on a minimal set of generators of $S_i$).

Therefore, in order to have the smallest Arf semigroup containing $S$, we must have $\textrm{Arf}(S) \in \sigma(E)$ where $E=\left\{ M_1,\ldots, M_n \right\}$ and $M_i$ is the multiplicity sequence associated to the Arf numerical semigroup $\textrm{Arf}(S_i)$.
\\
Now we need to find the matrix $$M(T)_{E}=\left(  \begin{matrix} 0 & p_{1,2} & p_{1,3} & \ldots & p_{1,n} \\ 
0 & 0 & p_{2,3} & \ldots & p_{2,n} \\ \ldots & \ldots & \ldots & \ldots & \ldots \\ 0 & 0 & 0 & \ldots& p_{n-1,n} \\ 0 & 0 & 0 & \ldots & 0 \\  \end{matrix}\right).$$that describes the tree of $\textrm{Arf}(S)$.

We recall that from the properties of good semigroups, it follows that there exists a minimal vector $\delta\in \mathbb{N}^n$ such that $\delta+\mathbb{N}^n \subseteq S$ (we will call this vector the conductor of $S$).

Suppose that $\delta=(c[1],\ldots,c[n])$. We denote by $$\textrm{Small}(S)=\left\{ \textbf{s}:  \textbf{0}<\textbf{s}\leq \delta \right\} \cap S,$$
the finite set of the small elements of $S$ (the elements of $S$ that are coordinatewise smaller  than the conductor). In \cite{DPMT} it is shown that $\textrm{Small}(S)$ describes completely the semigroup $S$ (in this paper we are not including the zero vector in $\textrm{Small}(S)$ to enlight the notations of the following procedures).

\begin{remark2}
We can recover the collection $E$ from $\textrm{Small}(S)$. In fact, the multiplicity sequence $M_i$ can be determined applying the Du Val algorithm to the set $ \left\{ s[i], s\in \textrm{Small}(S)\right\} \cup \left\{ c[i]+1 \right\} \subseteq S_i$. In order to apply the Du Val algorithm we  may have to add $c[i]+1$ because we can have $\gcd( \left\{ s[i], s\in \textrm{Small}(S)\right\}) \neq 1$. Because $c[i]$ and $c[i]+1$ belong to $S_i$, we know that $\textrm{Arf}(S_i)$ has conductor smaller than $c[i]$ and this implies that we only have to consider the elements that are smaller than $c[i]+1$.
\end{remark2}

Now, we notice that $p_{i,j} \leq \min(\textrm{pos}_E(i,c[i]),\textrm{pos}_{E}(j,c[j]))$ for all $i,j\in \left\{1,\ldots,n\right\}$, with $i<j$,  where we are using the notations of the previous section.
In fact,  if $\textrm{pos}_E(i,c[i]) \neq \textrm{pos}_E(j,c[j])$, we have already noticed that in an Arf semigroup containing $\delta$ the $i$-th and the $j$-th branches cannot be glued at a level greater than  $\min(\textrm{pos}_E(i,c[i]),\textrm{pos}_E(j,c[j]))$, then  $p_{i,j} \leq \min(\textrm{pos}_E(i,c[i]),\textrm{pos}_{E}(j,c[j]))$ .
If  $\textrm{pos}_E(i,c[i]) =\textrm{pos}_{E}(j,c[j])$ then we have $\delta_1=(c[1],\ldots,c[i]+1,c[i+1],\ldots,c[n]) \in S$, and $\textrm{pos}_E(i,c[i]+1)=\textrm{pos}_E(i,c[i])+1>\textrm{pos}_{E}(j,c[j])$.

 Therefore in  an Arf semigroup containing $\delta_1$ the $i$-th and the $j$-th branches cannot be glued at a level greater than  $$\min(\textrm{pos}_E(i,c[i])+1,\textrm{pos}_{E}(j,c[j]))=\textrm{pos}_{E}(j,c[j])=$$ $$=\min(\textrm{pos}_E(i,c[i]),\textrm{pos}_{E}(j,c[j]),$$ hence we have again $p_{i,j} \leq \min(\textrm{pos}_E(i,c[i]),\textrm{pos}_{E}(j,c[j]))$.\\
Furthermore, we always  have to take in account that $p_{i,j} \leq k_E(i,j)$ for all $i,j \in \left\{1,\ldots,n\right\}$.
\\ Now let us consider the following subset of $\left\{1,\ldots,n\right\}^2$, $$U_E(\textrm{Small}(S))=\left\{ (i,j): \textrm{pos}_E(i,s[i]) =\textrm{pos}_{E}(j,s[j]) \textrm{ for all } s \in \textrm{Small}(S) \right\}.$$
If $(i,j) \in \left\{1,\ldots,n\right\}^2 \setminus U_E(\textrm{Small}(S))$, and $i<j$ we can consider the following integers
$$ \textrm{MIN}_E(i,j,\textrm{Small}(S))=\min\left(k_E(i,j),\min\left\{ \min(\textrm{pos}_{E}(i,s[i]),\textrm{pos}_{E}(j,s[j])): s \in \textrm{Small}(S), \right. \right. $$ $$ \left. \left. \textrm{pos}_E(i,s[i]) \neq \textrm{pos}_{E}(j,s[j])  \right\}\right). $$
Notice that we need only to consider the vectors of $\textrm{Small}(S)$ because if $ s\in S $ then $ s_1=\min(s,\delta) \in \textrm{Small}(S)$ and we clearly have
$$ \min(\textrm{pos}_{E}(i,s[i]),\textrm{pos}_{E}(j,s[j])) \geq \min(\textrm{pos}_E(i,s_1[i]),\textrm{pos}_{E}(j,s_1[j])), $$
therefore $s_1 \in  \textrm{Small}(S)$ gives us more precise information on the ramification level than $s$ (it can happen that $\textrm{pos}_E(i,s_1[i])=\textrm{pos}_{E}(j,s_1[j])$ and $\textrm{pos}_{E}(i,s[i]) \neq \textrm{pos}_{E}(j,s[j]) $,  but only when
$ \min(\textrm{pos}_{E}(i,s[i]),\textrm{pos}_{E}(j,s[j])) \geq \min(\textrm{pos}_E(i,c[i]),\textrm{pos}_{E}(j,c[j]))$). 

Thus, if $T^1$ is an Arf semigroup of $\sigma(E)$ containing $S$, represented by $$M(T^1)_{E}=\left(  \begin{matrix} 0 & a_{1,2} & a_{1,3} & \ldots & a_{1,n} \\ 
0 & 0 & a_{2,3} & \ldots & a_{2,n} \\ \ldots & \ldots & \ldots & \ldots & \ldots \\ 0 & 0 & 0 & \ldots& a_{n-1,n} \\ 0 & 0 & 0 & \ldots & 0 \\  \end{matrix}\right) $$ we have: 
\begin{itemize}
\item  $a_{i,j} \leq \textrm{MIN}_E(i,j,\textrm{Small}(S))$ for $(i,j) \in  \left\{1,\ldots,n\right\}^2 \setminus U_E(\textrm{Small}(S));$

\item $a_{i,j} \leq \min(k_E(i,j),\textrm{pos}_{E}(i,c[i]))$, for $i \in U_E(\textrm{Small}(S))$ (we have $\textrm{pos}_E(i,c[i]) =\textrm{pos}_{E}(j,c[j])$).
\end{itemize}

Then we can finally deduce that the $p_{i,j}$ that defines the matrix of $\textrm{Arf}(S)$ are such that 
\begin{itemize}
\item  $p_{i,j} =\textrm{MIN}_E(i,j,\textrm{Small}(S))$, for $(i,j) \in  \left\{1,\ldots,n\right\}^2 \setminus U_E(\textrm{Small}(S));$

\item $p_{i,j} = \min(k_E(i,j),\textrm{pos}_{E}(i,c[i]))$, for $i \in U_E(\textrm{Small}(S))$ (we have $\textrm{pos}_E(i,c[i]) =\textrm{pos}_{E}(j,c[j])$),
\end{itemize}
and it is easy to see that the $p_{i,j}$ fulfill the condition of  Remark \ref{remark}.
\begin{remark23}
In other words we showed that $\textrm{Arf}(S)$ can be computed by computing $\textrm{Arf}(G)$ where:
$$ G=\textrm{Small}(S) \bigcup \left\{(c[1]+1,\ldots,c[i],c[i+1],\ldots,c[n]),\ldots,(c[1],\ldots,c[i]+1,c[i+1],\ldots, c[n]),\ldots,\right.$$ $$ \left.(c[1],\ldots,c[i],c[i+1],\ldots,c[n]+1)     \right\}.$$ 
\end{remark23}
\begin{example9}
	Let us consider the good semigroup $S$ with the following set of small elements,
	$$ \textrm{Small}(S)= \left\{ [ 5, 6, 5 ], [ 5, 10, 5 ], [ 5, 12, 5 ], [ 8, 6, 8 ], [ 8, 10, 8 ], [ 8, 12, 8 ], [ 8, 6, 10 ], [ 8, 10, 10 ], \right.$$ $$ \left.
	[ 8, 12, 10 ], [ 10, 6, 8 ], [ 10, 10, 8 ], [ 10, 12, 8 ], [ 10, 6, 10 ], [ 10, 10, 10 ], [ 10, 12, 10 ] \right\}.$$ 
	The conductor is $\delta=[10,12,10]$.
	First of all we need to recover from $\textrm{Small}(S)$ the collection of multiplicity sequences $E$.
	We have to apply the Du Val algorithm to the following sets:
	$$ \left\{ 5, 8, 10, 11\right\},  \left\{ 6, 10, 12, 13\right\} \textrm{ and } \left\{ 5, 8, 10, 11\right\}, $$
	 therefore we find that $E=\left\{ [5,3,2,1,1], [6,4,2,1,1], [5,3,2,1,1]  \right\}.$
	 
	 We have
	 $$ \textrm{pos}(\textrm{Small}(S))=\left\{\textrm{pos}_E(s): s \in \textrm{Small}(S)  \right\}=\left\{ [ 1, 1, 1 ], [ 1,2, 1 ], [ 1, 3, 1 ], [ 2, 1, 2 ], [ 2, 2, 2 ], \right.$$ $$ \left.[ 2, 3, 2 ], [ 2, 1,3 ], [ 2, 2, 3 ], 
	 [ 2, 3, 3 ], [ 3, 1, 2 ], [ 3, 2, 2 ], [ 3, 3, 2 ], [3, 1,3 ], [3, 2,3 ], [ 3,3, 3 ] \right\}.$$	
	 It is easy to check that $U_E(\textrm{Small}(S))=\emptyset.$ Thus we have
	 \begin{itemize}
	 	\item $p_{1,2}=\textrm{MIN}_E(1,2,\textrm{Small}(S))=\min(k_E(1,2)=4,1)=1,$ because we have the element $[1,2,1] \in \textrm{pos}(\textrm{Small}(S))$ corresponding to $s=[5,10,5] \in \textrm{Small}(S)$ such that $1=\textrm{pos}_E(1,s[1]) \neq \textrm{pos}_{E}(2,s[2])=2 $ and $\min(\textrm{pos}_{E}(1,s[1]),\textrm{pos}_{E}(2,s[2]))=1 $ .
	 	\item $p_{2,3}=\textrm{MIN}_E(2,3,\textrm{Small}(S))=\min(k_E(2,3)=4,1)=1,$ because we have the element $[1,2,1] \in \textrm{pos}(\textrm{Small}(S))$ corresponding to $s=[5,10,5] \in \textrm{Small}(S)$ such that $2=\textrm{pos}_E(2,s[2]) \neq \textrm{pos}_{E}(3,s[3])=1 $ and $\min(\textrm{pos}_{E}(2,s[2]),\textrm{pos}_{E}(3,s[3]))=1 $ .
	 	\item $p_{1,3}=\textrm{MIN}_E(1,3,\textrm{Small}(S))=\min(k_E(1,3)=+\infty,2)=2,$ because we have the element $[2,2,3] \in \textrm{pos}(\textrm{Small}(S))$ corresponding to $s=[8,10,10] \in \textrm{Small}(S)$ such that $2=\textrm{pos}_E(1,s[1]) \neq \textrm{pos}_{E}(3,s[3])=3 $ and $\min(\textrm{pos}_{E}(1,s[1]),\textrm{pos}_{E}(3,s[3]))=2 $, and we cannot find a smaller discrepancy.
	 	
	 	So the Arf closure of $S$ is described by the matrix $$ M(T)_E=\left( \begin{matrix} 0 &1 &2\\ 0& 0 &1 \\ 0 &0 &0 \\ \end{matrix}\right). $$ with $$E=\left\{M_1=[5,3,2,1,1],  M_2=[6,4,2,1,1] \textrm{ and } M_3=[5,3,2,1,1] \right\}.$$
	 \end{itemize}
 The following procedure, implemented in GAP, has as argument the set of small elements of a good semigroup and give as a result the Arf Closure of the given good semigroup. The Arf clousure is described by a list $[E,M(T)_E]$. 
 \begin{verbatim}
	gap> S:=[[5,6,5],[5,10,5],[5,12,5],[8,6,8],[8,10,8],[8,12,8],
	[8,6,10],[8,10,10],[8,12,10],[10,6,8],[10,10,8],[10,12,8],
	[10,6,10],[10,10,10],	[10,12,10]]; 
	
	[ [ 5, 6, 5 ], [ 5, 10, 5 ], [ 5, 12, 5 ], [ 8, 6, 8 ],
	 [ 8, 10, 8 ], [ 8, 12, 8 ], [ 8, 6, 10 ], [ 8, 10, 10 ],
	 	[ 8, 12, 10 ], [ 10, 6, 8 ], [ 10, 10, 8 ], [ 10, 12, 8 ],
	 	 [ 10, 6, 10 ], [ 10, 10, 10 ], [ 10, 12, 10 ] ]
	
	gap> ArfClosureOfGoodsemigroup(S); 
	 [ [ [ 5, 3, 2 ], [ 6, 4, 2 ], [ 5, 3, 2 ] ],
	  [ [ 0, 1, 2 ], [ 0, 0, 1 ], [ 0, 0, 0 ] ] ]
	\end{verbatim}
\end{example9}
\section{Bounds on the minimal number of vectors determining a given Arf semigroup } \label{section5}

Suppose that  $E$ is a collection of $n$ multiplicity sequences. Let $T\in \tau(E)$ and given a semigroup $S(T)$ in $\sigma(E)$,  we want to study the properties that a set of vectors $G(T) \subseteq \mathbb{N}^n$  has to satisfy to have  $S(T)=\textrm{Arf}(G(T))$, with the notations given in the previous section.  We call such a $G(T)$ a set of generators for $S(T)$. In particular we want to find bounds on the cardinality of a minimal set of generators for a $S(T) \in \sigma(E)$.

Since we want to find a $G(T)$ such that $\textrm{Arf}(G(T))$ is well defined, it has to satisfy the following properties:
\begin{itemize}

 \item For all $i=1,\ldots,n$
$$ \gcd(v[i]; v \in G(T))=1,$$
where $v[i]$ is the $i-$th coordinate of the vector $v \in G(T)$.
\item For all $i,j\in \left\{1,\ldots,n\right\}$, with $i<j$ there exists $v \in G(T)$ such that $v[i] \neq v[j]$. 
\end{itemize}
We recall that, given a Arf numerical  semigroup $S$, there is a uniquely determined smallest semigroup $N$ such that the Arf closure of $N$ is $S$. The minimal system of generators for such $N$ is called the Arf system of generators for $S$, or the set of characters of $S$. 

Now we want that  $\textrm{Arf}(G(T))$ is an element of $\sigma(E)$.
This implies that, when we apply the algorithm of Du Val to $G(T)[i]$, we have to find the $i$-th multiplicity sequence of $E$. This means that, if we call $S_i$  the Arf numerical semigroup corresponding to the projection on the $i$-th coordinate, we must have $G(T)[i] \subseteq S_i$ and furthermore $G(T)[i]$ has to contain a minimal system of generators for $S_i$.
 In fact, in \cite{Arf} it is proved that if we have $G=\left\{ g_1,\ldots,g_m \right\} \subseteq \mathbb{N}$ with $\gcd(G)=1$ then $G$ must contain the set of characters of the Arf closure of the semigroup $N=\langle G \rangle$.\\
So we need to recall a way to compute the characters of an Arf numerical semigroup. 
    
We suppose that $E=\left\{ M(1),\ldots, M(n) \right\} $.
Given $$ M(i)=[m_{i,1},\ldots, m_{i,M}],$$
we define the restricion number $r(m_{i,j})$ of $m_{i,j}$ as the number of sums $\displaystyle m_{i,q}=\sum_{h=1}^k {m_{i,q+h}}$ where $m_{i,j}$ appears as a summand.
With this notation we have that the characters of the multiplicity sequence $M(i)$ are the elements of the set (cf. \cite[ Lemma 3.1]{BDF2}) $$\textrm{Char}_E(i)=\left\{\sum_{k=1}^j{m_{i,k}} :r(m_{i,j})<r(m_{i,j+1})  \right\}.$$ 
Notice that, from our assumptions on $M$, it follows  that the last two entries in each $M(i)$ are $1$, and it is easy to see how it guarantees that we cannot find characters in correspondence of indices greater than $M$. We define  $\textrm{PChar}_E(i)=\left\{ j: r(m_{i,j})<r(m_{i,j+1}) \right\}$.

Given the collection $E$, we denote by $$V_E(j_1,j_2,\ldots,j_n)=\left[ \sum_{k=1}^{j_1}{m_{1,k}}, \sum_{k=1}^{j_2}{m_{2,k}},\ldots, \sum_{k=1}^{j_n}{m_{n,k}}\right].$$
Now, the elements of $G(T)$ must be of the type $V_E(j_1,j_2,\ldots,j_n)$ (in fact we noticed that when we project on the $k$-th coordinate we must find an element of the corresponding numerical semigroup that has the previous representation for some $j_k$).

From the previous remarks and notations we have the following property:
$$G(T)=\left\{ \textrm{Gen}(1)=V_E(j_{1,1},\ldots,j_{1,n}), \ldots, \textrm{Gen}(t)=V_E(j_{t,1},\ldots,j_{t,n})\right\}$$ are generators of a semigroup of $\sigma(E)$ if and only if  $$\textrm{PChar}_E(i) \subseteq \left\{ j_{1,i},\ldots,j_{t,i} \right\} \textrm{ for all }i=1,\ldots,n.$$
In particular we have found a lower bound for the cardinality of a minimal set of generators for a $S(T) \in \sigma(E)$.
In fact $G(T)$ has at least $C_E=\max\left\{ |\textrm{PChar}_E(i)|, i=1,\ldots,n\right\}$ elements.

Now we want to determine the generators of a given semigroup $S(T) \in \sigma({E})$. We have the following theorem.

\begin{theo3}
\label{T5}
Let $S(T) \in \sigma(E)$ with $$M(T)_{E}=\left(  \begin{matrix} 0 & p_{1,2} & p_{1,3} & \ldots & p_{1,n} \\ 
0 & 0 & p_{2,3} & \ldots & p_{2,n} \\ \ldots & \ldots & \ldots & \ldots & \ldots \\ 0 & 0 & 0 & \ldots& p_{n-1,n} \\ 0 & 0 & 0 & \ldots & 0 \\  \end{matrix}\right).$$ Denote by $P=\left\{ (q,r): p_{q,r}=k_E(q,r) \right\}$.
Then  $G(T)=\left\{\textrm{Gen}(1),\ldots,\textrm{Gen}(t) \right\} \subseteq \mathbb{N}^n$  is such that $\textrm{Arf}(G(T))=S(T)$  if and only if the following conditons hold
\begin{itemize}
\item $\textrm{Gen}(1)=V_E(j_{1,1},\ldots,j_{1,n}), \ldots, \textrm{Gen}(t)=V_E(j_{t,1},\ldots,j_{t,n})$ for some values of the indices $j_{1,1},\ldots,j_{t,n}$;
\item $\textrm{PChar}_E(i) \subseteq \left\{ j_{1,i},\ldots,j_{t,i} \right\}$ for all $i=1,\ldots,n$.

Furthermore, if we consider the following integer 
$$ \textrm{MIN}_{G(T)}(q,r)=\min\left( k_E(q,r),\min\left\{ \min(j_{p,q},j_{p,r}) : j_{p,q} \neq j_{p,r}, p=1,\ldots,t \right\} \right),$$
for the $(q,r)\in \left\{1,\ldots,n \right\}^2$ with $q<r$ and where it is well defined, we have:
\item for $(q,r) \in P$ we have either $ j_{p,q}=j_{p,r} $ for all $p=1,\ldots,t$, or  $ \textrm{MIN}_{G(T)}(q,r)$ is well defined and it equals $k_E(q,r)$;
\item $ \textrm{MIN}_{G(T)}(q,r)$ is well defined and it equals $p_{q,r}$, for all $(q,r) \notin P$.
\end{itemize}
\end{theo3}

\noindent \textbf{Proof} ($\Leftarrow$)
Suppose that we have $G(T)=\left\{\textrm{Gen}(1),\ldots,\textrm{Gen}(t) \right\} \subseteq \mathbb{N}^n$ satisfying the conditions of the theorem. The first two conditions ensure that if we apply the algorithm defined in the previous section on $G(T)$  it will produce an element of $\sigma(E)$.

Now it is easy, using the notations of  Theorem \ref{T45}, to show that $j_{p,q}=\textrm{pos}_E(q,\textrm{Gen}(p)[q])$ and from this it follows that, when  $\textrm{MIN}_{G(T)}(q,r)$ is well defined, it is equal to $\textrm{MIN}_E(q,r,G(T))$. Furthermore we have $U_E(G(T)) \subseteq P$. In fact we have $$ U_E(G(T))=\left\{(q,r)\in \left\{1,\ldots,n\right\}^2: \textrm{pos}_E(q,\textrm{Gen}(p)[q])=\textrm{pos}_E(r,\textrm{Gen}(p)[r]) \right.$$ $$ \left. \textrm{ for all } p=1,\ldots,t \right\}=\left\{(q,r)\in \left\{1,\ldots,n\right\}^2: j_{p,q}=j_{p,r} \textrm{ for all } p=1,\ldots,t \right\},$$
therefore if $(q,r) \in U_E(G(T))$ then $(q,r) \in P$, since $G(T)$ satisfy the fourth condition in the statement of the theorem (we cannot have $(q,r) \notin P$ because in this case $\textrm{MIN}_{G(T)}(q,r)=\textrm{MIN}_E(q,r,G(T))$ has to be well defined).
 So it will follows that, if $S(T^1)$ is $\textrm{Arf}(G(T))$ then $$M(T^1)_{E}=\left(  \begin{matrix} 0 & a_{1,2} & a_{1,3} & \ldots & a_{1,n} \\ 
0 & 0 & a_{2,3} & \ldots & a_{2,n} \\ \ldots & \ldots & \ldots & \ldots & \ldots \\ 0 & 0 & 0 & \ldots& a_{n-1,n} \\ 0 & 0 & 0 & \ldots & 0 \\  \end{matrix}\right) $$   where 
\begin{itemize}
\item $a_{i,j}=\textrm{MIN}_E(i,j,G(T))$ if $(i,j) \notin U_E(G(T))$;
\item $a_{i,j}=k_E(i,j)$ if $(i,j) \in U_E(G(T))$. \end{itemize}
Therefore if $(i,j) \notin P$ then $(i,j) \notin  U_E(G(T))$ and we have  $a_{i,j}=\min(\textrm{MIN}_E(i,j,G(T)))=\textrm{MIN}_{G(T)}(i,j)=p_{i,j}$.
If $(i,j) \in P$ then 
\begin{itemize}
\item if $(i,j) \in  U_E(G(T))$ then $a_{i,j}=k_E(i,j)$;
\item if $(i,j) \notin  U_E(G(T))$ then $a_{i,j}=\textrm{MIN}_E(i,j,G(T))=\textrm{MIN}_{G(T)}(i,j)=k_E(i,j)$, for the properties of the set $G(T)$ ($(i,j)\in P$).
\end{itemize}
So we showed that $\textrm{Arf}(G(T))=S(T)$.
Thus the proof of this implication is complete.

($\Rightarrow$)
It follows immediately by contradiction,  using the first part of the proof. \qed

\begin{example7} \label{ex7}

Suppose that we have $E= \left\{ M(1),M(2),M(3) \right\}$, where
$$ M(1)=[5,4,1,1],  M(2)=[6,4,1,1],  M(3)=[2,2,1,1]. $$
We have, $k_E(1,2)=3,  k_E(2,3)=2$ and $k_E(1,3)=2$. 

 We can define:
$$R(i)=[r(m_{i,1}),r(m_{i,2}),\ldots,r(m_{i,N}) ].$$
Notice that $r(m_{i,1})=0, r(m_{i,2})=1$. The values of $\textrm{PChar}(i)$ are the indices where this sequence has an increase (it can be easily shown that when the sequence has an increase we have $r(m_{i,j})=r(m_{i,j+1})-1$ cf. \cite[ Lemma 3.2]{BDF2}).
Furthermore $R(1)=[0,1,2,2,2,2]$, $R(2)=[0,1,2,3,2,2]$ and $R(3)=[0,1,1,2]$.
So $\textrm{PChar}_E(1)=\left\{1,2 \right\}, \textrm{PChar}_E(2)=\left\{1,2,3\right\}$ and $\textrm{PChar}_E(3)=\left\{1,3\right\}$.

Suppose that we want to find generators for the untwisted tree $T^1$ such that $ T^1_{E}=(2,1) $. We need at least three vectors because $C_E=3$. Consider the vectors  $\textrm{Gen}(1)=V_E(1,1,3), \textrm{Gen}(2)=V_E(2,3,2)$  and $\textrm{Gen}(3)=V_E(2,2,1)$.
The second condition, that guarantees that we have a tree belonging to $\tau(E)$, is satisfied. Furthermore $\textrm{MIN}_{G(T)}(1,2)=\min(k_E(1,2),2)=2$ ,  $\textrm{MIN}_{G(T)}(2,3)=\min(k_E(2,3),1)=1$, and $\textrm{MIN}_{G(T)}(1,3)=\min(k_E(1,3),1)=1=\min(d_1,d_2)$  where $G(T)=\left\{ \textrm{Gen}(1),\textrm{Gen}(2),\textrm{Gen}(3)\right\}$.
Thus we have $\textrm{Arf}(G(T))=S(T^1)$.
They are the vectors  $ \textrm{Gen}(1)=[5,6,5], \textrm{Gen}(2)=[9,11,4], \textrm{Gen}(3)=[9,10,2] $ which appeared in the Example \ref{ex6}.
\end{example7}
Now, we want to find an upper bound for the cardinality of a minimal set $G(T)$ such that $\textrm{Arf}(G(T)) \in \sigma(E)$. 
\begin{remark24} \label{remark24}
Suppose that $T^1$ is a twisted tree of $\tau(E)$, where $E$ is a collection of $n$ multiplicity sequences. Then there exists a permutation $\delta\in S^n$ such that $\delta(T^1)$ is an untwisted tree of $\tau(\delta(E))$. Then if $G$ is a set of generators for $\delta(T^1)$ then it is clear that we have
$$ \delta^{-1}(G)=\left\{ \delta^{-1}(g); g \in G \right\},$$
is a set of generators for the twisted tree $T^1$.
\end{remark24}
From the previous remark it follows that we can focus only on the untwisted trees of $\tau(E)$ to find an upper bound for the cardinality of $G(T)$. 

Our problem is clearly linked to the following question.
\begin{question}
Let us consider a vector $\textbf{d}=[d_1,\ldots,d_n] \in \mathbb{N}^n$. For all the $G \subseteq \mathbb{N}^{n+1}$  we denote by $\textrm{MIN}(G,i,j)$ the integers (if they are well defined)
$$ \textrm{MIN}(G,i,j)= \min\left\{ \min(g[i],g[j]): g\in G \textrm { with } g[i]\neq g[j] \right\}, $$
for all the $i < j$ and $i,j \in \left\{1,\ldots,n+1\right\}.$

We define a solution for the vector $\textbf{d}$ as a set $G \subseteq \mathbb{N}^{n+1}$ such that: 

$$\textrm{MIN}(G,i,j)=\min\left\{ d_i,\ldots,d_{j-1}\right\} \textrm{ for all } i<j. $$

Consider $n\in \mathbb{N}$ with $n \geq 1$ . Find the smallest  $t \in \mathbb{N}$, such that for all $[d_1,\ldots,d_{n}] \in \mathbb{N}^{n}$ there exists a solution with $t$ vectors. We denote such a number $t$ by $\textrm{NS}(n)$.
\end{question}
\begin{theo2} \label{thm}
Consider $n\in \mathbb{N}$ with $n \geq 1$ . Then $\textrm{NS}(n)=\left \lceil{ \log_2{(n+1)}}\right \rceil $, where $\left \lceil{d}\right \rceil=\min\left\{m\in \mathbb{N}: m\geq d \right\}$.
\end{theo2}
\noindent \textbf{ Proof} First of all we show that given an arbitrary vector $\textbf{d}$ of $\mathbb{N}^n$ we are able to find a solution of $\textbf{d}$ consisting of  $N=\left \lceil{ \log_2{(n+1)}}\right \rceil$ vectors.

We will do it by induction on $n$. The base of induction is trivial. In fact if $n=1$ then for each vector $[d_1]$ we find the solution $G=\left\{ [d_1,d_1+1]\right\}$ that has cardinality $\left \lceil{ \log_2{(1+1)}}\right \rceil=1$. Thus we suppose that the theorem is true for all the $m<n$ and we prove it for $n$. Let $\textbf{d}$ an arbitrary vector of $\mathbb{N}^n$. We fix some notations. Given a vector $\textbf{d}$, we will denote by $\textrm{sol}(\textbf{d})$ a solution with $\left \lceil{ \log_2{(n+1)}}\right \rceil$ vectors. We denote by $\textrm{Inf}(\textbf{d})=\min\left\{ d_i: i=1,\ldots,n \right\}$ and by $\textrm{Pinf}(\textbf{d})=\left\{i \in \left\{1,\ldots,n\right\}: d_i=\textrm{Inf}(\textbf{d}) 	\right\}$. We have $1\leq |\textrm{Pinf}(\textbf{d})|=k(\textbf{d}) \leq n$.

Suppose that $\textrm{Pinf}(\textbf{d})=\left\{ i_1 < i_2 < \dots < i_{k(\textbf{d})}\right\}$.
Then we can split the vector $\textbf{d}$ in the following $k(\textbf{d})+1$ subvectors:

$$ \begin{cases}\textbf{d}_1=\textbf{d}(1 \ldots i_1-1), \\ \textbf{d}_{j}=\textbf{d}(i_{j-1}+1 \ldots i_{j}-1) \textrm{ for } j=2,\ldots,k(\textbf{d}),\\ \textbf{d}_{k(\textbf{d})+1}=\textbf{d}(i_{k(\textbf{d})}+1\ldots n), \end{cases}$$
where with  $\textbf{d}(a \ldots b)$ we mean 
\begin{itemize}
\item $ \emptyset$ if $b<a$;
\item The subvector of $\textbf{d}$ with components between $a$ and $b$ if $a\leq b$.
\end{itemize}
Then the subvectors $\textbf{d}_j$ are either empty or with all the components greater than $\textrm{Inf}(\textbf{d})$.
We briefly illustrate with an example the construction of the subvectors $\textbf{d}_j$.
\begin{example4}
Suppose that $\textbf{d}=[2,3,2,2,5,4,5]$. Then $\textrm{Inf}(\textbf{d})=2$,  $\textrm{Pinf}(\textbf{d})=\left\{ 1,3,4\right\}$ and then  we have the four subvectors:
\begin{itemize}
\item $ \textbf{d}_1=\textbf{d}(1\ldots 0)=\emptyset,$
\item $ \textbf{d}_2=\textbf{d}(2\ldots 2)=[3]$,
\item $ \textbf{d}_3=\textbf{d}(4\ldots 3)=\emptyset $,
\item $ \textbf{d}_4=\textbf{d}(5\ldots 7)=[5,4,5]$.
\end{itemize}
\end{example4}
Then we can consider the list of $k(\textbf{d})+1$ subvectors:
$$ p(\textbf{d})=[\textbf{d}_1,\ldots,\textbf{d}_{k(\textbf{d})+1}],$$
and, because all the $\textbf{d}_i$ have length strictly less than $n$ we can find a solution for each of them with  $N=\left \lceil{ \log_2{(n+1)}}\right \rceil$ or less vectors. For the $\textbf{d}_i=\emptyset$ we will set $\textrm{sol}(\emptyset)=\left\{ [x]\right\} $, where $x$ is an arbitrary integer that is strictly greater than all the entries of $\textbf{d}$.
It is quite easy to check that the following equality holds:
\begin{equation} \label{form}n=k(\textbf{d})+\sum_{i=1}^{k(\textbf{d})+1}{\textrm{Length}(\textbf{d}_i}). \end{equation}
We associate to the list of vectors $p(\textbf{d})$ another list of vector $c(\textbf{d})$ such that 
$$ c(\textbf{d})=[\textbf{c}_1,\ldots,\textbf{c}_{k(\textbf{d})+1}],$$ 
where $\textrm{Length}(\textbf{c}_j)=\textrm{Length}(\textbf{d}_j)+1$ and the  entries of $\textbf{c}_j$ are all equal to $\textrm{Inf}(\textbf{d})$ for all $j=1,\ldots, k(\textbf{d})+1. $

Now we consider the set $ I(N)=\left\{0,1\right\}^N$. For each $\textbf{t} \in I(N)$ we denote by $O(\textbf{t})$ the number of one that appear in $\textbf{t}$. Because we have  $N=\left \lceil{ \log_2{(n+1)}}\right \rceil$ , it follows $$k(\textbf{d})+1 \leq n+1 \leq 2^N=|I(N)|,$$
therefore it is always possible to associate to each subvectors of the list $p(\textbf{d})$ distinct elements of $I(N)$.
We actually want to show that it is possible to associate to all the subvectors $\textbf{d}_i$ distinct vectors of $\textbf{t} \in I(N)$ such that $O(\textbf{t})\geq |\textrm{sol}(\textbf{d}_i)|$ (for $\textbf{d}_i=\emptyset$ we can also associate the zero vector).
We already know for the inductive step that all the $\textbf{d}_i$ have solutions with at most $N$ vectors. Suppose therefore that $m \leq N$.

 It is easy to see that
$$ |\left\{ \textbf{t} \in I(N) : O(t)\geq m \right\}|= \sum_{k=m}^N{\binom{N}{k}}.$$
Then we suppose by contradiction that in $p(\textbf{d})$ we have  $\sum_{k=m}^N{\binom{N}{k}}+1$ subvectors with solution with cardinality $m$. From the inductive step it follows that all these subvectors have at least length $2^{m-1}$, and from the formula \ref{form} it follows:

$$ n \geq \sum_{k=m}^N{\binom{N}{k}}+\left(\sum_{k=m}^N{\binom{N}{k}}+1\right)2^{m-1}  \Rightarrow n+1 \geq \left(\sum_{k=m}^N{\binom{N}{k}}+1\right)(1+2^{m-1}). $$
But we also have that:

$$\sum_{k=m}^N{\binom{N}{k}}+1 \geq 2^{N-m+1},$$
in fact $\sum_{k=m}^N{\binom{N}{k}}$ is the number of ways to select a subset of $\left\{ 1,\ldots,N\right\}$ of at least $m$ elements while there are  $2^{N-m+1}-1$ ways to select a subset  which contains at least $m$ elements and contains $\left\{ 1,2,\ldots,m-1\right\}$.

Therefore we can continue the inequality:

$$ n +1 \geq  2^{N-m+1} (1+2^{m-1})=2^N+2^{N-m+1}>2^N.$$
But $N=\left \lceil{ \log_2{(n+1)}}\right \rceil$  and therefore  $n+1 \leq 2^N$ and we find a contradiction.
Then in $\left\{ \textbf{t} \in I(N) : O(t)\geq m \right\} $ we have enough vectors to cover all the subvectors with solution with cardinality $m$. We still also have to exclude the following possibility.
Suppose that we have $x$ subvectors with solutions of cardinality $m_1$ and $y$  subvectors with solutions of cardinality $m_2 >m_1$. If $|\left\{ \textbf{t} \in I(N) : O(t)\geq m_1 \right\} |-x<y$ then it would not be possible to associate to all the subvectors of the second type an element $\textbf{t}$ of $I(N)$ with $O(\textbf{t}) \geq m_2.$ 
Indeed if this happen we would have:
$$ n \geq x+y-1+x\cdot2^{m_1-1}+y\cdot 2^{m_2-1} >x+y-1+(x+y)2^{m_1-1} \Rightarrow$$
$$ \Rightarrow n+1 \geq (x+y)(1+2^{m_1-1}) \geq \left(\sum_{k=m_1}^N{\binom{N}{k}}+1\right)(1+2^{m_1-1}), $$
and we already have seen that this is not possible.

Therefore we showed that we can consider a matrix $A$ with $N$ rows and $ k(\textbf{d})+1$ distinct columns with only zeroes and ones as entries and such that the $i$-th column of $A$ is a vector $\textbf{t}$ of $I(N)$ such that $O(\textbf{t}) \geq | \textrm{sol}(\textbf{d}_i)|$ for each $1\leq i \leq  k(\textbf{d})+1$.

Now we can complete the construction of a solution for $\textbf{d }$.
We consider a matrix $B$ with $N$ rows and $ k(\textbf{d})+1$  columns. We fill the matrix $B$ following these rules:
\begin{itemize}
\item If $A[i,j]=0$ then in $B[i,j]$ we put the vector $\textbf{c}_j$;
\item If $A[i,j]=1$ then in $B[i,j]$ we put an element of $ \textrm{sol}(\textbf{d}_j)$;
\item All the elements of  $\textrm{sol}(\textbf{d}_j)$ have to appear in the $j$-th column for all $j=1,\ldots, k(\textbf{d})+1$.
\end{itemize}
Then if we glue all the vectors appearing in the $i$-th row of $B$ for each $i=1,\ldots,N$ we obtain a solution $G$ for the vector $\textbf{d}$.
In fact if we consider $i_1,j_1$ such that $i_1<j_1$ we have two possibilities:
\begin{itemize}
\item $i_1$ and $j_1$ both correspond to elements in the $j$-th column of $B$. Then because  in this column we have either vectors of a solution for $\textbf{d}_j$  or costant vectors, it follows that they fulfill our conditions.
\item $i_1$ and $j_1$ correspond to elements in distinct columns. This implies that we must have $\textrm{MIN}(G,i_1,j_1)=\textrm{Inf}(\textbf{d})$. In fact, for construction, between two distinct subvectors we have an element equal to $\textrm{Inf}(\textbf{d})$ in $\textbf{d}$ forcing $\textrm{MIN}(G,i_1,j_1)=\textrm{Inf}(\textbf{d})$. Now suppose that $i_1$ and $j_1$ correspond respectively  to  elements in the $i$-th and $j$-th columns of $B$. Because we suppose $i \neq j$  we have that the $i$-th column  and the $j$-th column of  the matrix $A$ are distinct so there exists a $k$ such that $A[k,i]=0$ and $A[k,j]=1$ (or viceversa). This implies that in $B$ we have a row where in the $i$-th column there is the constant vector equal to $\textrm{Inf}(\textbf{d})$ while in the $j$-th column we have a vector corresponding to a solution of a subvectors of $\textbf{d}$ (that has all the components greater than $\textrm{Inf}(\textbf{d})$ by  construction). This easily implies that $\textrm{MIN}(G,i_1,j_1)=\textrm{Inf}(\textbf{d})$.
\end{itemize}
\begin{example5} Suppose that $\textbf{d}=[2,3,2,2,5,4,5]$. We have $n=7$, then we want to show that there exists a solution with three vectors. We have already seen that in this case we have:
$$ p(\textbf{d})=[\emptyset,[3],\emptyset,[5,4,5]].$$
We need to compute a solution for each entry of $p(\textbf{d})$. We have:
\begin{itemize}
\item $\textrm{sol}(\emptyset)=\left\{ [6]\right\}$ ( $6$ is greater than all the entries of $\textbf{d}$);
\item $\textrm{sol}([3])=\left\{ [3,4]\right\}$;
\item Let us compute a solution for $\textbf{f}=[5,4,5]$ with the same techniques. Because $\textrm{Length}(\textbf{f})=3$ we expect to find a solution with at most two vectors. We have:
$$ p(\textbf{f})=[[5],[5]],$$
and we have $\textrm{sol}([5])=\left\{ [5,6] \right\}$. Then in $I(2)$ we want to find two distinct vectors with at least an entry equal to one. We can choose $[1,1]$ and $[0,1].$ Therefore we have:
\[ A=\left( \begin{array}{cc}
1 & 0  \\
1 & 1  \\
\end{array} \right)  \textrm{ and } B=\left( 
 \begin{matrix}
  [5,6] & [4,4] \\
  [5,6] & [5,6] \\
 \end{matrix}
 \right). \]
Then  $\textrm{sol}([5,4,5])=\left\{ [5,6,4,4],[5,6,5,6]\right\}$.
\end{itemize}
Now we want to find in $I(3)$ four vectors $\textbf{t}_i$ for $i=1,\ldots,4$. We have free choice for the $\textbf{t}_1$ and $\textbf{t}_3$ , while we need $O(\textbf{t}_2) \geq 1$ and  $O(\textbf{t}_4) \geq 2$. For instance we choose $\textbf{t}_1=[0,0,0],\textbf{t}_2=[1,0,0],\textbf{t}_3=[1,1,0],\textbf{t}_4=[1,0,1]$. Then we have:

\[ A=\left( \begin{matrix}
0 & 1 & 1 & 1 \\
0 & 0 &  1 & 0\\
0 & 0 &0 &1 \\
\end{matrix} \right)  \textrm{ and } B=\left( 
 \begin{matrix}
 [2] & [3,4] & [6] & [5,6,4,4] \\
[2] & [2,2] &  [6] & [2,2,2,2]\\
[2] & [2,2] & [2] & [5,6,5,6]\\
 \end{matrix}
 \right). \]
Then a solution for $\textbf{d}$ is the set $$ G=\left\{[2, 3, 4, 6, 5, 6, 4, 4], [2, 2, 2, 6, 2, 2, 2, 2], [2, 2, 2, 2, 5, 6, 5, 6]\right\}.$$

\end{example5}
So we proved that $\textrm{NS}(n) \leq \left \lceil{ \log_2{(n+1)}}\right \rceil $. To prove that the equality holds we notice that for each $n$ a constant vector needs exactly $\left \lceil{ \log_2{(n+1)}}\right \rceil $ vectors in its solutions. \qed

Now  we can return to the  problem of determining an upper bound for the cardinality of $G(T)$.  We need another lemma:
\begin{Lemma2} \label{Lemma2}
Let $E=\left\{M(1),M(2) \right\}$ be a collection of two multiplicity sequences. Then, with the previous notations we have:
$$ k_E(1,2) \leq \min\left\{j:   j\in (\textrm{PChar}_E(1)\cup  \textrm{PChar}_E(2))\setminus (\textrm{PChar}_E(1)\cap  \textrm{PChar}_E(2)) \right\}.   $$
\end{Lemma2}
\noindent \textbf{Proof}
Let us choose an arbitrary element $t \in (\textrm{PChar}_E(1)\cup  \textrm{PChar}_E(2))\setminus (\textrm{PChar}_E(1)\cap  \textrm{PChar}_E(2))$. 
We want to show that $k_E(1,2) \leq t$. Suppose by contradiction that $t<k_E(1,2)$.
Without loss of generality we suppose that $t \in \textrm{PChar}_E(1)$. It follows that  $t \notin \textrm{PChar}_E(2)$ and we have:
$$ r(m_{1,t}) <r(m_{1,t+1}) \textrm{ and } r(m_{2,t}) \geq r(m_{2,t+1}). $$
Notice that if an entry of $M(1)$  has $m_{1,t+1}$ as a summand and it is not $m_{1,t}$, it is  forced to have $m_{1,t}$ as a summand too. So from  $r(m_{1,t}) <r(m_{1,t+1}) $ we deduce that in $M(1)$ there are no entries involving only $m_{1,t}$.
Similarly from  $r(m_{2,t}) \geq r(m_{2,t+1})$ we deduce that in $M(2)$ we must have at least one entry $m_{2,s}$ which involves $m_{2,t}$ as a summand but not $m_{2,t+1}$. 

Namely \begin{equation} m_{2,s}= \sum_{k=s+1}^t{m_{2,k}}. \label{equation} \end{equation}
Now, we have assumed that $t<k_E(1,2)$ hence $t+1\leq k_E(1,2)$. This imply that the untwisted tree $T$ such that $T_E=(t+1)$ is well defined. In $T$ we have the following nodes:
$$ (m_{1,s},m_{2,s}), \ldots, (m_{1,t},m_{2,t}), (m_{1,t+1},m_{2,t+1}).$$
Then from (\ref{equation}) and from the fact that the two branches are still glued at level $t+1$ it must follow that
$$ \label{eq} m_{1,s}= \sum_{k=s+1}^t{m_{1,k}} $$
and we have still noticed how it contradicts the assumption $ r(m_{1,t}) <r(m_{1,t+1}) $. \qed
\\ \\
Now we can prove the following result:
\begin{prop2} \label{prop2}
Let $E$ be a collection of $n$ multiplicity sequences. Then, if $S(T) \in \sigma(E)$, there exists $G(T) \subseteq \mathbb{N}^n$ with $\textrm{Arf}(G(T))=S(T)$ and  $|G(T)|=C_E+\left \lceil{ \log_2{(n})}\right \rceil $.
\end{prop2}
\noindent \textbf{Proof}
For the Remark \ref{remark24} it suffices to prove the theorem only for the untwisted trees. Therefore we suppose that $T_E=(d_1,\ldots,d_{n-1})$. First of all we have to satisfy the condition on the characters to ensure that $\textrm{Arf}(G(T)) \in \sigma(E)$. From the Lemma (\ref{Lemma2}) it follows that we can use $C_E$ vectors to satisfy all the conditions.
To see it, let us fix some notations.

Denote by $ \tau(i)=| \textrm{PChar}_E(i)|$ for all $i=1,\ldots,n$. Therefore $C_E=\max\left\{\tau(i), i=1,\ldots,n\right\}$. Suppose that
$$ \textrm{PChar}_E(i)=\left\{a_{i,1}<\dots<a_{i,\tau(i)}\right\},$$
and we define $$  L=\max\left(\bigcup_{i=1}^{n}{ \textrm{PChar}_E(i)}\right)+1. $$
For all $i=1,\ldots,n$ we consider the vector $J(i)=[a_{i,1},\ldots,a_{i,\tau(i)},L,\ldots,L] \in \mathbb{N}^{C_E}.$
Thus we can use the following set of vectors to satisfy the condition on the characters,
$$ G=\textrm{Gen}(1)=V_E(j_{1,1},\ldots,j_{1,n}), \ldots, \textrm{Gen}(C_E)=V_E(j_{C_E,1},\ldots,j_{C_E,n}),$$
where $j_{p,q}=J(q)[p]$ for all $p=1,\ldots,C_E$ and $q=1,\ldots,n$.
Now it is clear that we have $\textrm{PChar}_E(i) \subseteq \left\{ j_{1,i},\ldots,j_{C_E,i} \right\}$ for all $i=1,\ldots,n$. 

We  also need  to show that this choice does not affect the condition on $(d_1,\ldots,d_{n-1})$.
We define $P=\left\{ (q,r) \in \left\{1,\ldots,n\right\}^2: j_{p,q}=j_{p,r} \textrm{ for all }  p=1,\ldots,C_E \right\}.$ Thus for each $(q,r) \in P$ the previous vectors are compatible with the conditions on the element $p_{q,r}$ of $M(T)_E$.

For each $ (q,r) \notin P $, we consider $$p(q,r)=\min\left\{ p: j_{p,q} \neq j_{p,r} \right\}.  $$
Now, because the entries of the vectors $J(q)$ are in an increasing order, it is clear that we have 
 $$ \textrm{MIN}_{G}(q,r)=\min \left(k_E(q,r), \min \left\{ \min(j_{p,q},j_{p,r}) : j_{p,q} \neq j_{p,r} \right\}\right)=$$$$=\min \left( k_E(q,r), \min(j_{p(q,r),q},j_{p(q,r),r})\right),\textrm{for all }(q,r)\notin P.  $$ 

Furthermore, for the particular choice of the vectors $\textrm{Gen}(i)$ and of $L$, it is clear that from  $ j_{p(q,r),q} \neq j_{p(q,r),r}$, it follows that 
$$  \min(j_{p(q,r),q},j_{p(q,r)r}) \in (\textrm{PChar}_E(q)\cup  \textrm{PChar}_E(r))\setminus (\textrm{PChar}_E(q)\cap  \textrm{PChar}_E(r)), $$
and from the Lemma (\ref{Lemma2}), we finally have  $$\min(j_{p(q,r),q},j_{p(q,r),r}) \geq k_E(q,r) \textrm{ for all } (q,r) \notin P, $$
so the vectors $\textrm{Gen}(i)$  are compatible with our tree.

 Now from the Theorem (\ref{thm}) it follows that we can use $\left \lceil{ \log_2{(n})}\right \rceil $ vectors to  have a solution for the vector $[d_1,\ldots,d_{n-1}]$.  Adding the vectors corresponding to this solution to the previous $C_E$  we obtain a set $G(T)$ such that $\textrm{Arf}(G(T))=S(T)$. \qed

Notice that the first $C_E$ vectors may satisfy some conditions on the $d_i$, therefore it is possible to find $G(T)$ with smaller cardinality than the previous upper bound.

\begin{remark4}
Let us consider the Arf semigroup of the Example \ref{ex7}.

 It was $T=T_E=(2,1)$, where $$E= \left\{M(1)=[5,4,1,1], M(2)=[6,4,1,1], M(3)=[2,2,1,1] \right\},$$ with $$\textrm{PChar}_E(1)=\left\{1,2 \right\}, \textrm{PChar}_E(2)=\left\{1,2,3\right\} \textrm{ and } \textrm{PChar}_E(3)=\left\{1,3\right\}.$$
We found $G=\left\{ V_E(1,1,3),V_E(2,3,2),V_E(2,2,1)\right\}$ as a set such that $\textrm{Arf}(G)=S(T)$, and it is also minimal because we have $|G|=C_E$ and we clearly cannot take off any vector from it.
Using the strategy of the previous corollary we would find the vectors:
$$ \textrm{Gen}(1)=V_E(1,1,1), \textrm{Gen}(2)=V_E(2,2,3) \textrm{ and } \textrm{Gen}(3)=V_E(4,3,4),$$
that satisfy the conditions on the characters ($L=4$).

 We have to add  vectors that correspond to a solution for the vector $[2,1]$. For istance it suffices to consider $ [3,2,1]$ and therefore we will  add the vector $\textrm{Gen}(4)=V_E(3,2,1)$.
Notice how the set $G'=\left\{ V_E(1,1,1),V_E(2,2,3),V_E(4,3,4),V_E(3,2,1)\right\}$, with $|G'|>|G|$, is still minimal because we cannot remove any vector from it without disrupting the condition on the tree. 
Therefore we can have minimal sets of generators with distinct cardinalities.
\end{remark4}
\begin{example8}
Let us consider $$E=\left\{ M_1=[4,4,1,1],M_2=[6,4,1,1],M_3=[2,2,1,1],M_4=[3,2,1,1] \right\}.$$ We want to find a set of generators for the twisted tree $T$ of $\tau(E)$ such that:
$$ M(T)_E=\left(  \begin{matrix} 0 &2& 1& 2 \\ 0& 0& 1& 3 \\ 0& 0 &0 &1 \\ 0& 0 &0 &0 \\\end{matrix}\right).$$
First of all we notice that it is well defined because it satisfies the conditions given by the Remark \ref{rem} and we have $$ k(1,2)=2, k(1,3)=4,k(1,4)=2, k(2,3)=2, k(2,4)=3 \textrm{ and } k(3,4)=2.$$
We consider the permutation $\delta=(3,4)$ of $S^4$.
Then $\delta(T)$ is an untwisted tree of $\tau(\delta(E))$ and it is described by the vector $T_{\delta(E)}=(2,3,1)$.
We have:
\begin{itemize}
\item $ \textrm{PChar}_{\delta(E)}(1)=\left\{1,3  \right\};$
\item  $ \textrm{PChar}_{\delta(E)}(2)=\left\{1,2,3  \right\};$
\item $ \textrm{PChar}_{\delta(E)}(3)=\left\{1,2  \right\};$
\item  $ \textrm{PChar}_{\delta(E)}(4)=\left\{1,3  \right\}.$
\end{itemize}
Then with the vectors $V_{\delta(E)}(1,1,1,1),V_{\delta(E)}(3,2,2,3),V_{\delta(E)}(4,3,4,4)$, we satisfy the condition on the characters. We need to add the vectors corresponding to a solution for $[2,3,1]$. It suffices to add $V_{\delta(E)}(2,4,3,1)$.
Then $$ G(T)=\left\{ [4,6,3,2],[9,10,5,5],[10,11,7,6],[8,12,6,2]\right\},$$
is a set of generators for $\delta(T)$. Because $\delta^{-1}=(3,4)$, we have that 
$$ \delta^{-1}( G(T))=\left\{ [4,6,2,3],[9,10,5,5],[10,11,6,7],[8,12,2,6]\right\}$$
is a set of generators for the twisted tree $T$.
\end{example8}
 \begin{acknowledgements}
    The author would like to thank Marco D'Anna for his helpful comments and suggestions. Special thanks to  Pedro Garc\'ia-S\'anchez for his careful reading of an earlier version of the paper and for many helpful hints regarding the implementation in GAP of the presented procedures.
\end{acknowledgements}
\begin{otherlanguage}{english}

\end{otherlanguage}

GIUSEPPE ZITO-Dipartimento di Matematica e Informatica-Universit\`a di Catania-Viale Andrea Doria, 6, I-95125 Catania- Italy.

E-mail address: giuseppezito@hotmail.it

\end{document}